\documentclass[letterpaper, 10pt, journal, onecolumn, final]{IEEEtran}

\IEEEoverridecommandlockouts

\usepackage{cite}
\usepackage{amsmath,amssymb,amsfonts}
\usepackage{graphicx}
\usepackage{algorithm,algorithmicx,algpseudocode}
\usepackage{tcolorbox}
\usepackage{booktabs} %
\usepackage{subcaption}
\usepackage{caption}
\usepackage{soul}

\title{Singular Perturbations for Control Systems Safety}

\author{Stefano Di Gregorio,
  Guido Carnevale, Giuseppe
  Notarstefano \thanks{The work
      was partially funded by the FISA-2023-00210 project APACHE - CUP
      J53C25000520001, by Next Generation EU, National Recovery and Resilience Plan (NRRP), Mission
      4, Component 2, Investment 3.3; CUP J33C24001490009, and by IMA
      S.p.A.}  \thanks{The authors are with the Department of
    Electrical, Electronic and Information Engineering, University of
    Bologna, 40136, Bologna, Italy (e-mail:
    stefano.digregorio@unibo.it, guido.carnevale@unibo.it,
    giuseppe.notarstefano@unibo.it).}%
}

\newcommand{\cS}{\mathcal{S}}

\newcommand{\cA}{\mathcal{A}}
\newcommand{\cQ}{\mathcal{Q}}

\newcommand{\cC}{\mathcal{C}}

\newcommand{\cU}{\mathcal{U}}
\newcommand{\cK}{\mathcal{K}}

\newcommand{\R}{\mathbb{R}}
\newcommand{\N}{\mathbb{N}}

\newcommand{\x}{x}

\newcommand{\z}{z}

\newcommand{\pr}{\epsilon}

\newcommand{\iter}{{t}}

\newcommand{\ud}{^}

\newcommand{\xt}{\x(\iter)}

\newcommand{\zt}{\z(\iter)}

\newcommand{\xtp}{\dot{\x}(t)}

\newcommand{\ztp}{\dot{\z}(\iter)}

\newcommand{\tz}{\tilde{\z}}

\newcommand{\tztp}{\dot{\tz}(\iter)}

\newcommand{\tzt}{\tz(\iter)}

\newcommand{\tzta}{\tz(\tau)}

\newcommand{\zeq}{\z\ud{\text{eq}}}

\newcommand{\norm}[1]{\left \|#1 \right \|}

\newcommand{\T}{^\top}

\newcommand{\map}[3]{#1: #2 \rightarrow #3}

\newcommand{\until}[1]{\{1,\ldots,#1\}}
\newcommand{\col}{\textsc{col}}

\newcommand{\lip}{L}

\newcommand\oprocendsymbol{\hbox{$\square$}}
\newcommand\oprocend{\relax\ifmmode\else\unskip\hfill\fi\oprocendsymbol}
\def\eqoprocend{\tag*{$\square$}}
\def\er/{Erd\H{o}s-R\'enyi}

\newcommand{\uu}{u}
\newcommand{\ut}{\uu(\iter)}

\newcommand{\dyn}{f}
\newcommand{\sdyn}{f} 
\newcommand{\fdyn}{g}

\newcommand{\rs}{_{\text{rs}}}
\newcommand{\bl}{_{\text{bl}}}

\newcommand{\nablahbound}{\lip_{\cbf}}
\newcommand{\fredbound}{\lip\rs}

\newcommand{\drs}{\sdyn\rs}
\newcommand{\dbl}{\fdyn\bl}

\newcommand{\cortwoset}{\mathcal{B}}

\def\PL/{Polyak-\L{}ojasiewicz}

\newcommand{\xddots}{%
	\raise 4pt \hbox {.}
	\mkern 6mu
	\raise 1pt \hbox {.}
	\mkern 6mu
	\raise -2pt \hbox {.}
}

\graphicspath{{figs/}}

\def\er/{Erd\H{o}s-R\'enyi}

\newcommand{\scale}{1}

\newcommand{\safeset}{\cS_x}
\newcommand{\safesetint}{\cS}

\newcommand{\cbf}{h}

\newcommand{\pol}{\pi} 
\newcommand{\lipp}{\lip_{\pol}}

\newcommand{\scl}{\tilde{\sdyn}_\pol}
\newcommand{\fcl}{\tilde{\fdyn}_\pol}
\newcommand{\zeqcl}{\zeq_{\pol}}

\newcommand{\out}{a}

\newcommand{\rob}{\eta}

\newcommand{\cCf}{\cC_{\text{\tiny full}}}

\newcommand{\lp}{\lambda_\pr}

\newcommand{\xipr}{\xi_{\pr}}

\newcommand{\res}{c}
\newcommand{\resd}{c^{\text{des}}}

\def\citeDL/{\cite[Proposition~6.1.2]{bertsekas2015convex}}
 
\newtheorem{theorem}{Theorem}

\newtheorem{definition}{Definition}
\newtheorem{lemma}{Lemma}

\newtheorem{assumption}{Assumption}

\begin{document}

\maketitle
\thispagestyle{empty}
\pagestyle{empty}

\begin{abstract}
  Control Barrier Functions (CBFs) offer a powerful theory for enforcing safety in dynamical systems.
  Yet, applying this framework to complex systems is challenging when subsystems are feedback-interconnected and operate at different timescales.
  Such structures arise naturally in multiphysics systems and in control architectures where physical plants interact with optimization, estimation, and control algorithms.
  In many of these settings, safety involves only part of the system, motivating the design of reduced-order certificates.
  Inspired by singular perturbation theory, we establish sufficient conditions, involving separation between slow and fast dynamics, under which a reduced-order safety certificate can be lifted to the full interconnected system.
  Specifically, we construct a composite CBF ensuring forward invariance of a lifted subset of the original safe set. We also show that the certified safe set expands monotonically as the tuning parameter accelerating the fast dynamics is reduced.
  Numerical results validate the approach on two robotic systems: a planar 2-DoF manipulator with motor dynamics and a 6-DoF UR10e manipulator controlled by a suboptimal optimization algorithm.
\end{abstract}

\section{Introduction}
The growing complexity of autonomous systems operating in real-world environments has led to extensive research in the field of safety-critical control, that is, the design of control strategies that guarantee the satisfaction of safety specifications, as reported in the recent survey \cite{hsu2023safety}.
Among the many approaches proposed in the literature, Control Barrier Functions (CBFs) provide a systematic means of enforcing safety in nonlinear systems by imposing suitable conditions on their dynamics (see, e.g., the works~\cite{ames2019control,xiao2023safe,wang2025safe}).
Their applicability extends to a wide range of domains, including robotic systems~\cite{ferraguti2022safety}, autonomous vehicles~\cite{alan2023control}, and human-robot interaction scenarios \cite{landi2019safety}.
Beyond enforcing individual safety specifications, several works have investigated the interaction between safety requirements and other control objectives, such as stability.
In particular, the combination of Control Lyapunov Functions (CLFs) and CBFs has been explored through optimization-based formulations~\cite{ames2016control} and analytical constructions for safe stabilization~\cite{ong2019universal}. 
Related works address robustness to disturbances through the input-to-state safety (ISSf) framework \cite{kolathaya2018input}, as well as difficulties arising with simultaneous stabilization and constraint satisfaction, including topological deadlocks via hybrid CLF-CBF control~\cite{matias2025hybrid} and parametric model uncertainty through composite adaptive certificates~\cite{kamaldar2026composite}.
\\
As modern control architectures increasingly involve multiple coupled components, constructing safety certificates for the overall system has become a critical challenge (see the early references~\cite{sloth2012compositional} and~\cite{sloth2012existence}).
To address safety in interconnected systems, several compositional approaches have been investigated. 
In particular, \cite{lyu2022small} establishes small-gain conditions to verify ISSf properties of interconnected systems.
The small-gain approach has been extended in~\cite{lyu2023small} to deal with high-relative-degree constraints.
This direction is further developed in~\cite{long2024sufficient}, which provides sufficient and necessary conditions for compositional ISSf of nonlinear interconnected systems.
Similarly, \cite{nejati2022compositional} constructs barrier certificates for interconnected stochastic hybrid systems by combining subsystem-level certificates through small-gain conditions.
Other approaches combine multiple safety constraints through Boolean compositions, see~\cite{molnar2023composing} and~\cite{yang2025control}, including nested and combinatorial specifications~\cite{ong2025combinatorial}.

A complementary research direction focuses on exploiting reduced-order models to simplify the synthesis and implementation of safety-critical controllers while retaining formal guarantees for the overall system (see, e.g., the recent survey~\cite{cohen2024safety}).
In~\cite{molnar2023safety}, input-to-state safe backup sets transfer reduced-order safety guarantees to full-order systems under input constraints.
Several contributions exploit structural properties of the dynamics to establish safety guarantees.
For systems admitting a strict-feedback representation, the paper~\cite{taylor2022safe} develops a backstepping-based (see, e.g.,~\cite{freeman1993backstepping}) construction extending safety conditions from a virtual controller to the actual input.
Related approaches exploit the mechanical structure of robotic systems through energy-based constructions~\cite{singletary2021safety}, while other techniques investigate model-free tracking of safe reference velocities~\cite{molnar2021model}. 
More recently, \cite{cohen2025safety} introduced a framework based on projection mappings and simulation functions to establish safety relationships between a reduced-order model and the corresponding full system.
From a complementary perspective, the authors of~\cite{de2026lifting} investigate the existence of safety certificates for full-order systems starting from safe reduced-order models. 
Related methods for constructing safety certificates in the presence of high-dimensional dynamics include the use of High-Order CBFs (HOCBFs)~\cite{xiao2021high}, which address the challenge of high relative degree by recursively constructing safety conditions through higher-order derivatives of the safety function, partial input-output feedback linearization-based constructions~\cite{cohen2024constructive} and rectified CBFs~\cite{ong2024rectified}.
\\
  The aforementioned approaches provide valuable tools for safety certification by exploiting specific structural properties of the system dynamics. A complementary research direction concerns the safety of interconnected systems evolving at different time scales, where the interactions between subsystems must be explicitly accounted for.

Classical singular 
perturbation (or timescale separation) theory (see, e.g., the pioneering work~\cite{kokotovic1968singular} and the surveys~\cite{kokotovic1976singular,kokotovic1999singular,abdelgalil2023multi}) provides a powerful framework to analyze the relationship between the reduced-order and full-order dynamics from the point of view of stability and convergence properties.
Timescale separation (or singular perturbation) theory provides a suitable framework for analyzing not only systems that couple different physical domains (e.g., the electrical and mechanical components of robotic systems) but also physical plants interconnected with online control and/or optimization algorithms.
The work in \cite{allibhoy2023control} develops CBF-based reformulations of gradient flows for constrained nonlinear optimization, ensuring forward invariance of the feasible set.
The papers~\cite{labar2019constrained,williams2022practically,williams2023semi,williams2024semiglobal,williams2024experimental,wang2026logarithmic,ghaffari2026safe} rely on stability arguments grounded in singular perturbation and averaging methods to guarantee \emph{practical safety} when the invariance certificate is dynamically and approximately imposed with the well-known extremum-seeking mechanism (see, e.g.,~\cite{krstic2000stability}).
Moreover, the works~\cite{delimpaltadakis2025feedback,delimpaltadakis2026safe} propose a safe feedback optimization method 
that enforces state constraints at all times, see also a related application to optimal power flow pursuit~\cite{colot2024optimal}.
More closely related to the present work, \cite{dall2026local} investigates locally implemented safety filters for networked systems through a two-timescale design and analysis.
See also the conference paper~\cite{krishnamurthy2026control}, which investigates safety guarantees for singularly perturbed systems by using an approach based on input-to-state stability.

The main contribution of this paper is a singular perturbation
framework for certifying safety of nonlinear feedback-interconnected
systems.
Specifically, we establish sufficient conditions under which a CBF designed for a reduced-order model can be
lifted to get a safety certificate for the full interconnected system.
As customary in singular perturbation approaches, we distinguish
between \emph{fast} and \emph{slow} dynamics. The fast dynamics is
assumed to possess a tunable parameter to modulate its speed and a
stable equilibrium manifold, while the slow one collapses to the
reduced-order model when the fast state is on the equilibrium
manifold.
Under sufficient timescale separation, we characterize a safe invariant set for the interconnected system controlled by a safe policy designed on the reduced-order model.
Drawing on composite Lyapunov methods for stability analysis, we construct a family of composite CBFs, depending on the tuning parameter. We combine the reduced-order safety certificate with a specially designed shaping function obtained from a Lyapunov function associated with the fast system.
The resulting forward invariant set quantifies the admissible deviations of the fast states from their
equilibrium manifold.
We further establish a nesting property of these safe sets as timescale separation increases and characterize the exact recovery of the reduced-order safety certificate in the singular limit.
The framework accommodates nonlinear feedback interconnections without requiring a cascade structure or inversion-based feedback laws.
The safety controller can therefore be implemented using the reduced-order model, while the resulting safety guarantees hold true for the full interconnection.

A preliminary version of this work appeared in the short conference paper~\cite{di2026safe}, which established the basic safety-lifting result.
The present paper extends that result by introducing a broader family of composite CBFs based on nonlinear transformations of the boundary-layer Lyapunov function and by establishing the nesting property of the corresponding certified safe sets.
It also provides self-contained proofs of the theoretical results and a more detailed geometric and analytical interpretation of the safety guarantees.
Finally, the numerical validation is expanded to include suboptimal safe control of a robotic manipulator simulated in a high-fidelity physics engine.

The paper is organized as follows. 
In Section~\ref{sec:preliminaries}, we review some preliminaries on safe control for continuous-time systems and the related safety guarantees.
In Section~\ref{sec:framework}, we introduce the framework under consideration and our goal.
In Section~\ref{sec:examples}, we present motivating examples to illustrate the application of the proposed framework.
In Section~\ref{sec:result}, we establish the main theoretical results of the paper, accompanied by a graphical and analytical interpretation based on a toy example, and we formally prove them in the Appendix.
Finally, in Section~\ref{sec:simulations}, we evaluate the proposed approach with numerical simulations and virtual experiments.

\paragraph*{Notation}
$\R_{+}$ is the set $\{x \! \in \! \R \! \mid \! x\ge 0\}$.
We denote by $\col(v_1,\dots,v_N)$ the vertical concatenation of the vectors $v_1, \dots, v_N$.
The symbol $\norm{\cdot}$ denotes the Euclidean norm for vectors and its induced norm for matrices.
Let $X \subseteq \R^{n}$ and $x \in \R^n$. Then $\|x\|_X := \inf_{y \in X}\norm{x - y}$.
A continuous function $\map{\alpha}{\R_+}{\R_+}$ is a class-$\cK_\infty$ function, denoted by $\alpha \in \cK_\infty$, if it is strictly increasing, $\alpha(0) = 0$, and $\lim_{r \to \infty} \alpha(r) = \infty$.
Given a continuously differentiable function $f:\R^n\times\R^m\to\R$, we denote by $\nabla_1 f(x,u) =
  \frac{\partial f(x,u)}{\partial x}
^\top \in \R^n$ and $\nabla_2 f(x,u) = 
  \frac{\partial f(x,u)}{\partial u}
^\top \in\R^m$ its gradients with respect to the first and second arguments, respectively.
For vector-valued functions, $\nabla$ denotes the corresponding Jacobian matrix.
Given a set $A$, we denote its closure by $\overline{A}$.

\section{Preliminaries on Safe Control}
\label{sec:preliminaries}

We start by recalling the theoretical foundations of Control Barrier Functions (CBFs) for continuous-time nonlinear systems of the form
\begin{align}\label{eq:generic_system}
  \xtp = \dyn\left(\xt,\ut\right),
\end{align}
where $\xt \in \R^{n}$ is the system state at time $\iter \in \R_+$, $\ut \in \cU \subseteq \R^{m}$ is the control input at time $\iter$, and $\dyn: \R^{n} \times \R^{m} \to \R^{n}$ describes the system dynamics.
The next definition introduces the notion of a zero-superlevel set of a function.
\begin{definition}[Zero-Superlevel Set]\label{def:zero_superlevel}
  Given a function $\cbf: \R^{n} \to \R$, the zero-superlevel set of $\cbf$ is defined as $\cC := \{\x \in \R^{n} \mid \cbf(\x) \geq 0\}$, while $\partial \cC := \{\x \in \R^{n} \mid \cbf(\x) = 0\}$.
  \oprocend
\end{definition}
Let us now 
introduce $\safeset \subseteq \R^n$ as the safe set of system~\eqref{eq:generic_system}, i.e., the set of states that satisfy the prescribed safety specifications. 
Accordingly, we say that system~\eqref{eq:generic_system} is safe with respect to $\safeset$ if there exists a set $\cC \subseteq \safeset$ that is forward invariant for system~\eqref{eq:generic_system}.
A popular tool for constructively certifying the safety of system~\eqref{eq:generic_system} is given by CBFs \cite{ames2016control, ames2019control, cohen2024safety}.
\begin{definition}[Control Barrier Function]
  \label{def:cbf}
  A continuously differentiable function $\cbf: \R^{n} \to \R$ with zero-superlevel set $\cC$ is said to be a \emph{Control Barrier Function (CBF)} for system~\eqref{eq:generic_system} on $\cC$ if there exists a class-$\cK_\infty$ function $\alpha \in \cK_\infty$ such that
  \begin{align*}
    \sup_{\uu \in \cU} \nabla \cbf(x)\T \dyn(\x,\uu) \geq -\alpha(\cbf(\x)),
  \end{align*}
  for all $\x \in \cC$. \oprocend
\end{definition}
The above condition characterizes the set of control inputs that are compatible with the prescribed safety constraint. As formalized by the following result, which directly follows from \cite[Corollary 1]{cohen2024safety}, any sufficiently regular feedback controller $\pi: \R^{n} \to \cU$ satisfying the CBF condition for all $\x \in \cC$ guarantees that the safe set $\cC$ is forward invariant for system~\eqref{eq:generic_system}, thus certifying its safety with respect to $\safeset$.
\begin{lemma}
  \label{lemma:cbf}
  Consider a continuously differentiable function $\cbf: \R^{n} \to \R$ with zero-superlevel set $\cC$. Assume that $\dyn$ is locally Lipschitz and that $\nabla \cbf(x) \neq 0$ for all $x \in \partial \cC$.
  If a locally Lipschitz feedback controller $\pol: \R^{n} \to \cU$ satisfies
  \begin{align*}
    \nabla \cbf(x)\T \dyn(\x,\pol(x)) \geq -\alpha(\cbf(\x)),
  \end{align*}
  for all $x \in \cC$, and some function $\alpha \in \cK_\infty$, then $\cC$ is forward invariant for the closed-loop system
   \begin{align*}
    \xtp = \dyn(\xt,\pi(\xt)). \eqoprocend
  \end{align*}
\end{lemma}

\section{Singular Perturbations for Safety: Problem Statement}
\label{sec:framework}

\begin{figure}[H]
  \centering
  \includegraphics[scale=\scale]{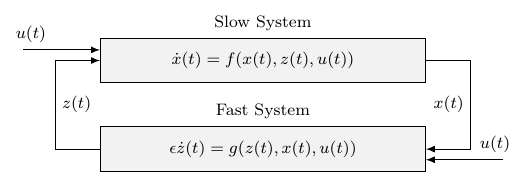}
  \caption{Block diagram representation of~\eqref{eq:interconnected_system_generic}.}
  \label{fig:scheme}
\end{figure}

In this section, we outline our procedure based on singular perturbations to lift reduced-order safety certificates in feedback-interconnected systems.
Indeed, we focus on systems of the form
\begin{subequations}\label{eq:interconnected_system_generic}
  \begin{align}
    \xtp &= \sdyn(\xt,\zt,\ut)\label{eq:slow_system_generic}
    \\
    \pr\ztp &= \fdyn(\zt,\xt,\ut),\label{eq:fast_system_generic}
  \end{align}
\end{subequations}	
where $\xt \in \R^{n}$ and $\zt \in \R^{p}$ are the system states at
time $\iter \in \R_+$, $\ut \in \cU \subseteq \R^{m}$ is the input applied at time
$\iter$, $\sdyn: \R^n \times \R^p \times \R^m \to \R^n$ and
$\fdyn: \R^p \times \R^n \times \R^m \to \R^p$ are the vector fields
of the slow and fast subsystems, respectively, and $\pr > 0$ is a tunable parameter to
modulate the relative speeds of these two subsystems.
By borrowing the typical nomenclature of timescale separation (or singular perturbations) for stability analysis (see, e.g., the recent survey~\cite{abdelgalil2023multi} or~\cite[Ch.~11]{khalil2002nonlinear}), we call $\xt$ and~\eqref{eq:slow_system_generic} the \emph{slow} state and the \emph{slow} system, respectively; $\zt$ and~\eqref{eq:fast_system_generic} are the \emph{fast} state and the \emph{fast} system, respectively.
Fig.~\ref{fig:scheme} describes system~\eqref{eq:interconnected_system_generic} in terms of block diagrams.

\subsection{Reduced-Order Safety Certificates}

\begin{figure}[H]
    \centering
    \includegraphics[width=0.67\textwidth, trim=2cm 0 3cm 1cm, clip]{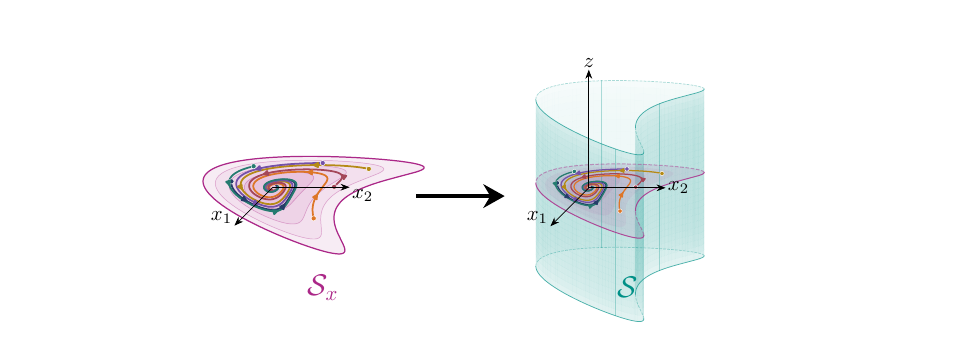}
    \caption{
        Graphical representation of the safe set $\safesetint$ (cf. \eqref{eq:safesetint}).}
    \label{fig:two_sets_concept}
  \end{figure}

\begin{figure}[H]
  \centering
  \includegraphics[scale=\scale]{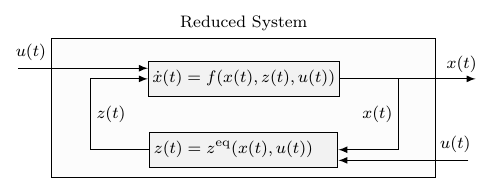}
  \caption{Block diagram representation of the reduced system~\eqref{eq:reduced}.}
  \label{fig:scheme_reduced}
\end{figure}

  Our approach builds upon a safety certificate established for a reduced-order representation of the interconnected dynamics.
  Its lifting requires three main ingredients: safety specifications involving only the slow state, an equilibrium manifold for the fast dynamics, and a feedback policy satisfying a robust safety condition for the reduced system. 
  \\
  First, we consider safety constraints that depend only on the slow state $\x$.
  Hence, we denote by $\safeset \subseteq \R^{n}$ the safe set of the slow subsystem~\eqref{eq:slow_system_generic} and define the safe set $\safesetint \subseteq \R^{n} \times \R^{p}$ for the overall system~\eqref{eq:interconnected_system_generic} as
  \begin{align}\label{eq:safesetint} 
    \safesetint := \{(x,z) \in \R^{n} \times \R^{p} \mid x \in \safeset\}.
  \end{align}
  Figure~\ref{fig:two_sets_concept} provides a graphical representation of the safe set $\safesetint$ in the full state space, obtained by extending the slow-state safe set $\safeset$ along the unconstrained fast-state directions.
  The second requirement concerns the existence of a parameter-dependent equilibrium function for the fast subsystem~\eqref{eq:fast_system_generic}.
  \begin{assumption}\label{ass:equilibrium}
  There exists a continuously differentiable function $\zeq: \safeset \times \cU \to \R^{p}$ such that
  \begin{align*}%
    0 &= \fdyn(\zeq(\x,u),\x,u),
  \end{align*}
  for every $\x \in \safeset$ and $u \in \cU$. \oprocend 
\end{assumption}
Under this assumption, the fast dynamics can be replaced by their equilibrium manifold to obtain a lower-dimensional description of the overall interconnection. 
Specifically, restricting the full system~\eqref{eq:interconnected_system_generic} to the manifold $\{(x,z,u) \in \safeset \times \R^{p} \times \cU \mid z = \zeq(x,u)\}$ yields the reduced-order vector field $\drs: \safeset \times \cU \to \R^n$ defined as
\begin{align}\label{eq:reduced}
  \drs(x,\uu) := \sdyn(x,\zeq(x,u),\uu).
\end{align}
This representation describes the evolution of the slow state $\x$ in the idealized regime where the fast state instantaneously tracks the equilibrium $\zeq(\x,\uu)$ associated with the current slow state $x$ and input $u$.
Hence, as illustrated in Fig.~\ref{fig:scheme_reduced}, the corresponding architecture is obtained by replacing the fast subsystem with the static map $\zeq$.

Finally, we require a feedback controller $\pol: \cC \to \cU$ that enforces a robust safety condition for the reduced-order model~\eqref{eq:reduced} on a compact zero-superlevel set $\cC \subseteq \safeset$.
\begin{assumption}\label{ass:cbf}
  There exist a continuously differentiable function
  $\cbf: \R^{n} \to \R$ with a compact zero-superlevel set
  $\cC \subseteq \safeset$, a continuously differentiable and $\lipp$-Lipschitz continuous function $\pol: \cC \to \cU$, an
  $\lip_{\alpha}$-Lipschitz continuous function
  $\alpha \in \cK_\infty$, and a safety margin $\rob > 0$ such
  that
  \begin{align}\label{eq:CBF}
    \nabla\cbf(x)\T \drs(x,\pol(x)) &\geq - \alpha(\cbf(x)) + \rob,
  \end{align}
  for all $x \in \cC$.\oprocend
\end{assumption}
Unlike the standard CBF condition (see Definition~\ref{def:cbf}), the inequality~\eqref{eq:CBF} incorporates a strictly positive safety margin $\rob$. 
This additional margin provides robustness against the discrepancy between the reduced-order dynamics $\xtp = \drs(x,\uu)$ and the actual dynamics~\eqref{eq:slow_system_generic} incorporating the fast state $z$.

\subsection{Boundary Layer System}

Following standard practice in singular perturbation theory (see, e.g.,~\cite{abdelgalil2023multi}), we complement the reduced-order description by introducing the so-called \emph{boundary layer} system.
  This auxiliary model captures the evolution of the fast state $\zt$ when the slow state $\xt$ is held constant, providing insights into the stability properties of the fast subsystem.
  Specifically, let $\bar{x} \in \safeset$ and $\bar{u} \in \cU$ be fixed values of the slow state and input, respectively.
  By holding these quantities fixed and introducing the stretched fast time variable $\tau := \iter/\pr$, we isolate the fast dynamics from the evolution of the slow subsystem.
  We then express the fast state relative to its equilibrium as $\tzta := \z(\tau) - \zeq(\bar{x},\bar{u}) \in \R^{p}$, yielding the boundary layer system
  \begin{equation}\label{eq:bl}
	\frac{\partial \tzta}{\partial \tau} = \dbl(\tzta,\bar{x},\bar{u}),
\end{equation}
with $\tzta \in \R^p$, where $\dbl: \R^p \times \safeset \times \cU \to \R^p$ is defined as
\begin{align}\label{eq:dbl_def}
  \dbl(\tz,\x,\uu) := \fdyn(\tz + \zeq(\x,\uu),\x,\uu).%
\end{align}

\begin{figure}[H]  
  \centering
  \includegraphics[scale=\scale]{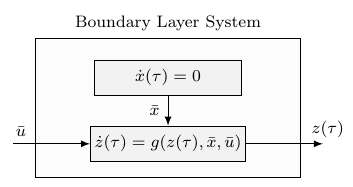}
  \caption{Block diagram representation of~\eqref{eq:bl} in the original coordinates $(x,z)$.}
  \label{fig:scheme_bl}
\end{figure}

\begin{figure}[htpb]
  \centering
  \includegraphics[width=0.67\textwidth, trim=0 0 0.5cm 1cm, clip]{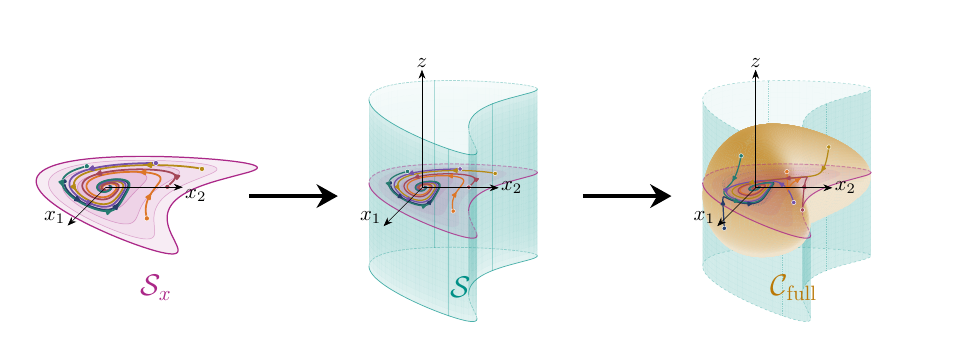}
  \caption{
      Conceptual illustration of the proposed framework. Starting from a safety certificate for the reduced-order dynamics, the goal is to construct a forward-invariant set $\cCf \subseteq \safesetint$ within the safe set $\safesetint$ (cf. \eqref{eq:safesetint}) of the full system accounting for the fast-state deviations from the equilibrium manifold $\zeq(\xt,\ut)$.
  }
  \label{fig:sets_three_concept}
\end{figure}

As depicted in Fig.~\ref{fig:scheme_bl}, the boundary-layer system~\eqref{eq:bl} represents system~\eqref{eq:interconnected_system_generic} in the limit $\pr \to 0$, which is equivalent to freezing the slow state and the input at a given $(\bar{x}, \bar{u}) \in \safeset \times \cU$.
We characterize the stability properties of system~\eqref{eq:bl} through the following Lyapunov conditions. 
\begin{assumption}\label{ass:bl}
  There exists a continuously differentiable function $U: \R^p \to \R$ such that
	\begin{subequations}\label{eq:U_generic}
		\begin{align}
			b_1 \norm{\tz}^2 \leq U(\tz) &\leq b_2\norm{\tz}^2
      \label{eq:U_first_bound_generic}
			\\
      \nabla U(\tz)\T \dbl(\tz,\bar{\x},\bar{\uu}) 
      &\leq 
      - b_3\norm{\tz}^2\label{eq:U_minus_generic}
			\\
      \norm{\nabla U(\tz)} &\leq b_4\norm{\tz},
      \label{eq:lipshcitz_hessian_u}
		\end{align}
	\end{subequations}
	for all $\tz \in \R^p$, $\bar{\x} \! \in \! \safeset$, $\bar{u} \in \cU$, and some $b_1, b_2, b_3, b_4 > 0$. \oprocend
\end{assumption}
Together, Assumptions~\ref{ass:cbf} and~\ref{ass:bl} characterize two complementary properties of the interconnected system~\eqref{eq:interconnected_system_generic}: (i) robust safety of the reduced-order dynamics and (ii) exponential stability of the equilibrium manifold for the boundary-layer dynamics.
The proposed approach combines these two properties to construct a safety certificate for the interconnected system within the full safe set $\safesetint$ defined in~\eqref{eq:safesetint}, as conceptually illustrated in Fig.~\ref{fig:sets_three_concept}.

We conclude this section by imposing some regularity conditions on the dynamics functions $\sdyn$ and $\fdyn$ of system~\eqref{eq:interconnected_system_generic}.
\begin{assumption}\label{ass:lip}
    The functions $\sdyn$ and $\fdyn$ are locally Lipschitz continuous. Moreover, $\sdyn(x,\cdot,\pol(x))$ is $L_{\sdyn}$-Lipschitz continuous uniformly with respect to $x \in \cC$. \oprocend
\end{assumption}

\section{Motivating Scenarios}\label{sec:examples}

To illustrate the applicability of our framework across diverse domains, this section presents representative scenarios characterized by distinct multi-timescale interconnections. 
Specifically, we focus on both physical plant-level couplings (such as electromechanical scenarios) and algorithmic interconnections, wherein physical systems interact with an online optimization and/or a control algorithm.

\subsection{Planar Robotic Arm with Actuator Dynamics}
\label{eq:motor_pendulum}

First, we focus on a fully actuated planar robotic arm with two degrees of freedom, in which the timescale separation is intrinsically determined by the nature of the mechanical and electrical domains. 
Specifically, the mechanical dynamics of the system are described by the standard rigid-body equations
\begin{align}
  \label{eq:pendulum}
  \begin{bmatrix}
    \dot{q} \\
    \dot{\omega}
  \end{bmatrix} =
  \begin{bmatrix}
      \omega \\
      M(q)^{-1} \Big( K_I I - B \omega - C(q, \omega)\omega - G(q) \Big)
    \end{bmatrix},  
\end{align}
where the terms $M(q) \in \R^{2 \times 2}$, $C(q, \omega) \in \R^{2 \times 2}$, and $G(q) \in \R^2$ represent the inertia matrix, the Coriolis matrix, and the gravity vector, respectively, and the diagonal and positive definite matrices $B, K_I \in \R^{2 \times 2}$ model the joint viscous friction and the motor torque constants, respectively, while
\begin{align*}
  q := \begin{bmatrix}\theta_1\\ \theta_2\end{bmatrix} \in \mathbb{T}^2, \quad 
  \omega := \begin{bmatrix}\omega_1\\ \omega_2\end{bmatrix} \in \R^2, \quad
  I := \begin{bmatrix}I_1\\ I_2\end{bmatrix} \in \R^2,
\end{align*}
are the vectors of angular positions, angular velocities, and motor currents, respectively. 
The mechanical plant is coupled with the electrical dynamics of the two motors, namely
\begin{align}
  \label{eq:motor}
    L \dot{I} = V - RI - K_\omega \omega,
  \end{align}
where $V = [V_1, V_2]^\top \in \R^2$ is the voltage applied to the motors, while $L$, $R$, $K_\omega \in \R^{2 \times 2}$ are positive definite diagonal matrices representing the inductance, resistance, and back-electromotive force constants, respectively.
By assuming that $L$ can be arbitrarily designed, we treat it as the tunable parameter $\pr$ in~\eqref{eq:interconnected_system_generic}.
Hence, we note that~\eqref{eq:pendulum} and~\eqref{eq:motor} can be framed within the interconnected system~\eqref{eq:interconnected_system_generic} with 
\begin{align*}
 x := \col(q, \omega), \quad \z := I, \quad \uu := V,
\end{align*}
and equilibrium function $I_\text{eq}: \mathbb{T}^2 \times \R^2 \times \R^2 \to \R^2$ defined as
\begin{align}
  I_\text{eq}(\col(q,\omega),V) := R^{-1} \left( V - K_\omega \omega \right).
\end{align}
Moreover, since $R \succ 0$, each equilibrium $I_\text{eq}(\col(q,\omega),V)$ is globally exponentially stable uniformly in $(\col(q,\omega), V)$.
Hence, Assumptions \ref{ass:equilibrium} and~\ref{ass:bl} are satisfied. 
Assumption~\ref{ass:cbf} thus corresponds to finding a safe controller for the reduced mechanical dynamics, namely 
\begin{align}\label{eq:motor_pendulum_red}
  \begin{bmatrix}
    \dot{q} 
    \\
    \dot{\omega}
  \end{bmatrix}
   =
  \begin{bmatrix}
    \omega \\
    M(q)^{-1} \Big( K_I R^{-1} V - B_{f} \omega - C(q, \omega)\omega - G(q) \Big)
  \end{bmatrix},
\end{align}
where $B_{f} := B + K_I R^{-1} K_\omega \in \R^{2 \times 2}$ is the effective damping matrix incorporating the back-electromotive force of the motors.

\subsection{Physical Dynamics with Safety Control Algorithm}
\label{sec:dynamic_safe_control}

Second, we consider the case in which the slow subsystem~\eqref{eq:slow_system_generic} is a given plant, while the fast one~\eqref{eq:fast_system_generic} is an optimization algorithm aimed at finding a safe control input that minimizes the distance from a desired but possibly unsafe control action. 
We formulate this scenario by introducing, for all $t \in \R_+$, an optimization problem of the form
\begin{subequations}\label{eq:problem}
  \begin{align}
    &\min_{u \in \cU} \quad \ell(u,\xt) 
    \\
    &\text{s.t. } \quad \nabla \cbf(\xt)\T \sdyn(\xt,u) \geq -\alpha(\cbf(\xt)) + \rob,  
  \end{align}
\end{subequations}
where $\ell: \R^m \times \R^n \to \R$ is the cost function
for a desired policy $\pi^{\text{des}}: \R^n \to \cU$, and $\cbf: \R^n \to \R$ is a CBF for the plant dynamics $\sdyn$.
In real applications, optimization problems like~\eqref{eq:problem} cannot be solved exactly in real time due to limited computational resources.
Instead, their solution estimates are iteratively refined by optimization algorithms (e.g., a gradient flow) that we can model as dynamical systems.
Thus, by generally denoting by $\cA: \R^{p} \times \R^{n} \to \R^{p}$ the specific algorithm flow used to dynamically address problem~\eqref{eq:problem}, the resulting interconnected system is given by 
\begin{subequations}\label{eq:plant_algo}
  \begin{align}
    \xtp &= \sdyn(\xt,\out(\zt))\label{eq:plant}
    \\
    \pr\ztp &= \cA(\zt,\xt),\label{eq:algo}
  \end{align}
\end{subequations}
where $\pr > 0$ modulates the speed of algorithm~\eqref{eq:algo}, while $\out: \R^p \to \cU$ is an output map selecting only the current solution estimate associated with the overall algorithmic state $\zt$, which may contain additional variables (e.g., multiplier estimates and momentum terms).
In this setup, Assumption~\ref{ass:cbf} is satisfied when problem~\eqref{eq:problem} admits a solution function $\uu_\star: \cC \to \cU$ that is continuously differentiable in the set $\cC$.
Assumption~\ref{ass:equilibrium} is satisfied provided that $\cA$ admits an equilibrium function $\zeq: \cC \to \R^p$ such that $\out(\zeq(x)) = \uu_\star(x)$ for all $x \in \cC$.
Assumption~\ref{ass:bl} is achieved with exponential stability of $\zeq(x)$ for system 
\begin{align*}
  \ztp = \cA(\zt,x),
\end{align*}
uniformly in $\x \in \cC$, namely, when the considered algorithm converges exponentially fast to the solution of~\eqref{eq:problem} for all $x \in \cC$.
Finally, Assumption~\ref{ass:lip} corresponds to Lipschitz continuity of the plant and algorithm dynamics $\sdyn$ and $\cA$.

\section{Singular Perturbations for Safety: Main Result}
\label{sec:result}

In this section, we establish the main safety result for the interconnected system~\eqref{eq:interconnected_system_generic} by combining the reduced-order safety guarantees encoded in $\cbf$ (cf. Assumption~\ref{ass:cbf}) and the stability properties of the boundary-layer system~\eqref{eq:bl} encoded in $U$ (cf. Assumption~\ref{ass:bl}).
\begin{theorem}\label{th:SP}
  Consider system \eqref{eq:interconnected_system_generic} and let Assumptions~\ref{ass:equilibrium},~\ref{ass:cbf},~\ref{ass:bl},~\ref{ass:lip} hold.
  Then, there exist $\bar{\pr} > 0$ and a family of composite barrier functions $V_{\pr}: \cC \times \R^p \to \R$ parameterized by $\pr > 0$, such that the following properties hold. 
  \begin{enumerate}
    \item \textbf{Forward invariance.} For all $\pr \in (0,\bar{\pr})$, the zero-superlevel set
    \begin{align}
      \cCf(\pr) := \left\{(x,z) \in \cC \times \R^p \mid V_{\pr}(x,z) \geq 0\right\}
      \label{eq:certified_safe_set}
    \end{align}
    is forward invariant for the closed-loop system
  \begin{align*}
    \xtp &= \sdyn(\xt,\zt,\pol(\xt))
    \\
    \pr\ztp &= \fdyn(\zt,\xt,\pol(\xt)).
  \end{align*}
  Equivalently, for all $(x(0),z(0)) \in \cCf(\pr)$, it holds $(x(t),z(t)) \in \cCf(\pr) \subseteq \safesetint$ for all $t \in \R_+$.
  \item \textbf{Nesting of certified safe sets.} The family of certified safe sets $\cCf(\pr)$, indexed by $\pr$, is monotonically non-increasing with respect to set inclusion, i.e.,
  \begin{align*}
    \cCf(\pr_2) \subseteq \cCf(\pr_1),
  \end{align*}
  for all $0 < \pr_1 \leq \pr_2 < \bar{\pr}$.
  \item \textbf{Recovery of the reduced-order certificate.} 
  For every compact subset $\cK$ of $\cC \times \R^p$, it holds
  \begin{align*}
    \lim_{\pr \to 0^+} \sup_{(x,z) \in \cK} | V_\pr(x,z) - \cbf(x) | = 0.
  \end{align*}

  Moreover, for every compact subset $\cK$ of the strict superlevel set $\{(x,z) \in \cC \times \R^p \mid h(x) > 0\}$, there exists $\bar{\pr}_{\cK} \in (0, \bar{\pr})$ such that $\cK \subseteq \cCf(\pr)$ for all $\pr \in (0, \bar{\pr}_\cK)$.

  Finally, the closure of the union of all certified safe sets recovers the entire region $\cC \times \R^p$, namely
  \begin{align*}
    \overline{\bigcup_{\pr \in (0, \bar{\pr})} \cCf(\pr)} = \cC \times \R^p.\eqoprocend
  \end{align*}
  \end{enumerate}
\end{theorem}
The proof of Theorem \ref{th:SP} is provided in the Appendix, with the three results established in Appendices \ref{sec:proof2}, \ref{sec:proof_cor1}, and \ref{sec:proof_cor2}, respectively.

The idea of the proof is to show that we are able to construct a composite CBF $V_\pr$ for the interconnected system~\eqref{eq:interconnected_system_generic}. 
The composite function combines the reduced-order CBF $\cbf$ (cf. Assumption \ref{ass:cbf}) with a suitable nonlinear transformation of the boundary-layer Lyapunov function $U$ (cf. Assumption~\ref{ass:bl}).
  In analogy with singular perturbation theory for stability, the first statement of Theorem~\ref{th:SP} establishes that, for sufficiently small values of $\pr$, the sets $\cCf(\pr)$ are forward invariant for the interconnected system \eqref{eq:interconnected_system_generic}, thus ensuring safety since $\cCf(\pr) \subseteq \safesetint$ by construction.
  Beyond forward invariance, the second statement characterizes how the certified safe region depends on timescale separation. 
  Specifically, the nesting property establishes that decreasing $\pr$ (i.e., increasing the timescale separation) progressively enlarges the set of initial conditions from which safety of the interconnected system can be guaranteed.
  Finally, the third statement characterizes the asymptotic behavior of the certified safe sets as $\pr \to 0^+$. 
  In particular, the composite CBF $V_\pr$ converges uniformly to the reduced-order CBF $\cbf$ on compact subsets of $\cC \times \R^p$.
  Consequently, every compact subset of the strict superlevel set $\{(x,z) \in \cC \times \R^p \mid h(x) > 0\}$ is contained in a forward-invariant set $\cCf(\pr)$ for sufficiently small $\pr > 0$.
  In this sense, the certified safe sets asymptotically recover the entire region $\cC \times \R^p$ as expressed by the closure of their union.

\subsection{Graphical Interpretation of Theorem \ref{th:SP}}

To illustrate the theoretical guarantees provided by the main results, we consider as a toy example the following two-dimensional singularly perturbed system
\begin{subequations}\label{eq:toy}
  \begin{align}
    \dot{x}(t) &= z(t) - x(t) - w(x(t)),\label{eq:slow_toy} 
    \\
    \epsilon \dot{z}(t) &= k\left(
      u(t) + x(t) + w(x(t)) - z(t)
    \right),\label{eq:fast_toy}
  \end{align}
\end{subequations}
with $x(t)$, $z(t)$, $\ut \in \R$, $w(x(t)) = 5 \sin(x(t))$ and $k = 5$. 
The fast subsystem~\eqref{eq:fast_toy} admits the equilibrium mapping $\zeq(x,u) = u + x + w(x)$, which yields the single-integrator reduced-order model $\dot{x}(t) = u(t)$.
We consider the physical safety specification for the slow state $x$ defined by the CBF $\cbf(x) = 1 - \frac{x^2}{\beta^2}$, which encodes the safe set
\begin{align*}
  \safeset = \left\{
    x \in \R \middle| 1 - \tfrac{x^2}{\beta^2} \geq 0
  \right\},
\end{align*}
with $\beta = \pi$. 

\begin{figure*}[t]
    \centering

    \begin{minipage}[t]{0.48\textwidth}
        \centering
        \includegraphics[width=\linewidth]{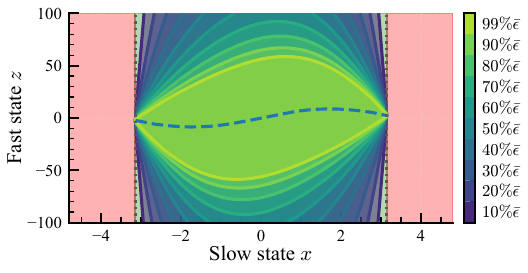}
        \label{fig:nested_sets}
    \end{minipage}
    \hfill
    \begin{minipage}[t]{0.48\textwidth}
        \centering
        \includegraphics[width=\linewidth]{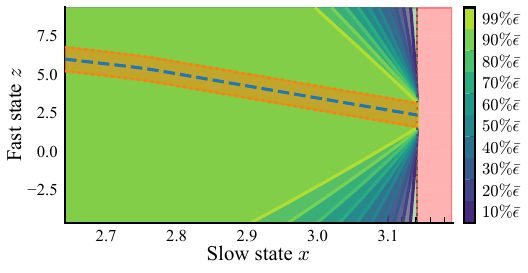}
        \label{fig:nested_sets_zoom}
    \end{minipage}

    \vspace{-1em}

    \begin{minipage}[t]{0.48\textwidth}
        \centering
        \includegraphics[width=\linewidth]{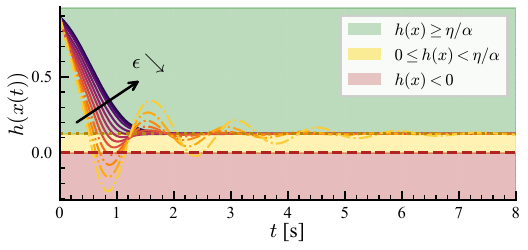}
        \label{fig:cbf_timeseries}
    \end{minipage}
    \hfill
    \begin{minipage}[t]{0.48\textwidth}
        \centering
        \includegraphics[width=\linewidth]{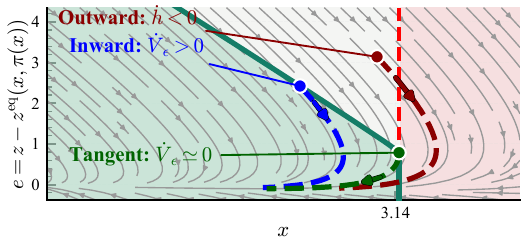}
        \label{fig:boundary_phase_portrait}
    \end{minipage}
    \vspace{-1em}
    \caption{
        Illustration of the certified sets and their associated
        dynamics. \textbf{Top left:} Sets $\cCf(\pr)$ for 10 values of
        $\pr\in(0,\bar{\pr})$, illustrating the nesting property
        established in Theorem~\ref{th:SP}.
        The blue dashed curve denotes the fast equilibrium manifold.
        \textbf{Top right:} Zoomed view highlighting the common region around
        the fast equilibrium manifold.
        \textbf{Bottom left:} Time evolution of $\cbf(x(t))$ for different values of $\pr$.
        \textbf{Bottom right:} Phase portrait illustrating outward-pointing dynamics
        near the physical safety boundary and inward-pointing or
        tangent dynamics at the boundary of the lifted certified set.}
    \label{fig:toy_example_summary}
\end{figure*}

The safe feedback policy $u = \pi(x)$ is obtained by smoothly approximating a CBF-QP that minimally modifies the nominal controller $\pi_{\text{des}}(x) = 2 \sin(x)$ to enforce the safety constraint~\eqref{eq:CBF} in Assumption \ref{ass:cbf} with $\rob = 0.5$.
For this system, the boundary-layer error dynamics are exponentially stable and admit the quadratic Lyapunov function $U(\tz) = (z - \zeq(x,\pi(x)))^2$.
The figure illustrates the geometry of the certified sets with the timescale parameter threshold $\bar{\pr} = 0.101$ computed by following the steps of the proof of Theorem~\ref{th:SP}.
In Figure~\ref{fig:toy_example_summary} we provide graphical representations that help visualize the concepts highlighted in the following discussion.

\paragraph{\textbf{Geometry and nesting of the certified sets}}
The top left panel of Figure \ref{fig:toy_example_summary} illustrates the certified safe sets $\cCf(\pr) = \{(x,z) \in \R^2 \mid V_{\pr}(x,z) \geq 0\}$ for different values of $\pr$, together with the nonlinear equilibrium manifold $z = \zeq(x,\pol(x))$. 
All certified sets lie within the physical safe region $\cS = \{(x,z) \in \R^n \times \R^p \mid |x| \leq \beta\}$.
As established by Theorem \ref{th:SP}, the certified safe sets are nested with respect to the singular perturbation parameter $\pr$, i.e., $\cCf(\pr_2) \subseteq \cCf(\pr_1)$ for all $0 < \pr_1 \leq \pr_2 < \bar{\pr}$, that is, decreasing $\pr$ enlarges the region of initial conditions for which safety of the full interconnected system can be guaranteed.
Specifically, the zoomed view in the top right panel highlights a common neighborhood of the equilibrium manifold in which, in light of the aforementioned construction, the composite barrier function $V_\pr$ coincides with the physical CBF $\cbf$. Consequently, within this neighborhood, the certified safe sets $\cCf(\pr)$ coincide with the physical safe set $\cS$.
\paragraph{\textbf{Influence of timescale separation on safety}} 
The bottom left panel of Figure \ref{fig:toy_example_summary} shows the evolution in time of $\cbf(x(t))$ for different values of $\pr$, starting from a common initial condition $(x(0),z(0)) = (1.0, \zeq(1.0) + 6.0)$. 
As the timescale parameter increases, the fast-state mismatch persists for longer and may produce larger deviations from the reduced-order dynamics. Thus, the safety condition may be violated when the fast subsystem is not sufficiently fast to compensate for the initial mismatch. This illustrates the importance of timescale separation in extending reduced-order safety guarantees to the full interconnected system. 
\paragraph{\textbf{Forward invariance through safety lifting}} The bottom right panel of Figure \ref{fig:toy_example_summary} shows the phase-plane representation in the $(x,e)$ coordinates, where $e = z - \zeq(x,\pol(x))$ is the fast-state mismatch. In these coordinates, the slow dynamics satisfy $\dot{x} = \pol(x) + e$, showing explicitly how the fast-state mismatch acts as a perturbation to the reduced-order dynamics. Although the controller $\pi(x)$ satisfies the reduced-order CBF condition, this additional term may cause the vector field to point outward at the physical safety boundary $\cbf(x) = 0$, as illustrated by the red trajectory.  
In contrast, the composite certificate $V_\pr$ accounts for the fast-state deviation and identifies a forward-invariant subset of the physical safe region. 
The blue and green trajectories show, respectively, inward-pointing and initially tangent vector fields at the boundary of the composite safe set (i.e., $V_\pr(x,z) = 0$). These observations provide a geometric interpretation of Theorem \ref{th:SP}: rather than requiring the entire physical safe set to be forward invariant, the lifting procedure constructs an explicitly characterized full-state region in which safety is guaranteed using a controller designed solely on the basis of the reduced-order model. 

\section{Numerical Simulations}\label{sec:simulations}

In this section, we showcase the applicability of the proposed framework in the scenarios introduced in Section~\ref{sec:examples}.
 Namely, we first consider a fully actuated 2-DoF robotic arm with joint motor dynamics, and then we consider a 6-DoF robotic manipulator controlled through an iterative optimization algorithm. 
These examples illustrate two distinct sources of timescale separation: the coupling between the mechanical and electrical dynamics of a robotic arm, and the cyber-physical interconnection between the physical system and a dynamic safety controller (i.e., an optimization algorithm).
In both cases, the safety controller is designed based solely on the reduced-order dynamics, while the simulations and the virtual experiments account for the full-order dynamics. 

\subsection{Planar Robotic Arm with Actuator Dynamics}\label{sec:simulations_robotic_arm}

First, we consider the fully actuated 2-DoF robotic arm described in Section~\ref{eq:motor_pendulum}, where the mechanical dynamics~\eqref{eq:pendulum} are coupled with the electrical dynamics~\eqref{eq:motor} of the joint motors.
In accordance with our proposed framework, the safety controller is designed based solely on \eqref{eq:motor_pendulum_red}.
The physical parameters of the robotic arm and the joint motors are provided in Table \ref{tab:robot_parameters}, assuming identical links and motors. 
\begin{table}[H]
\centering
\caption{Manipulator and Actuator Parameters}
\label{tab:robot_parameters}
\begin{tabular}{|c|c|c|c|c|c|c|c|c|}
\hline
\textbf{Symbol} & $m$ & $l$ & $J$ & $b$ & $R$ &$K_I$ & $K_\omega$ \\ \hline
\textbf{Value}  & 0.5 & 1.0 & 0.167 & 0.5 & 0.6 & 0.4 & 0.4        \\ \hline
\textbf{Units}  & kg & m & kg$\cdot$m$^2$ & $\frac{\mathrm{N}\cdot\mathrm{m}\cdot\mathrm{s}}{\mathrm{rad}}$ & $\Omega$ & $\frac{\mathrm{N}\cdot\mathrm{m}}{\mathrm{A}}$ & $\frac{\mathrm{V}\cdot\mathrm{s}}{\mathrm{rad}}$ \\ \hline
\end{tabular}
\end{table}
As discussed in Section~\ref{sec:examples}, the motor inductance $L$ is treated as a tunable parameter $\pr$ to modulate the timescale separation between the mechanical and electrical dynamics.
We consider $N_c=9$ safety constraints comprising both the joint limits and the obstacle avoidance requirements.
Specifically, the joint limits are
$q_{\min} = [-185^\circ, -80^\circ]^\top$ and $q_{\max} = [185^\circ, 80^\circ]^\top$ and the corresponding safety functions are defined as
\begin{align*}
  \bar{h}_1(q) &:= q_1 - q_{\min,1}
  ,\quad
  \bar{h}_2(q) := q_{\max,1} - q_1
  \\
  \bar{h}_3(q) &:= q_2 - q_{\min,2} 
  ,\quad
  \bar{h}_4(q) := q_{\max,2} - q_2.
\end{align*}
In addition, we consider $5$ circular obstacles with radius $r = 0.3$ m centered at the workspace coordinates $c_1 = [1.1, -1.7]^\top$, $c_2 = [2.2, 0.0]^\top$, $c_3 = [1.9, 0.6]^\top$, $c_4 = [1.3, 1.7]^\top$, and $c_5 = [-0.4, 2.1]^\top$.
The corresponding safety functions are
\begin{align*}
  \bar{h}_{4+k}(q) := \norm{p(q) - c_k}^2 - r^2, \quad k \in \left\{1, \dots, 5\right\},
\end{align*}
where $p(q) \in \R^2$ is the position of the end-effector.
Since these constraints have a relative degree of two with respect to the input, we construct a High-Order CBF (see, e.g., \cite{xiao2021high}) for each $j\in\{1,\dots,N_c\}$, namely, given $\gamma > 0$, we define
\begin{align}
  h_j(q, \omega) := \nabla \bar{h}_j(q)\T \omega + \gamma \bar{h}_j(q).
\end{align}
To simplify the filter design, all constraints are aggregated into a single CBF $h(q, \omega)$ using the Log-Sum-Exp formulation (see, e.g., \cite{safari2025safe}), namely
\begin{align}
  h(q, \omega) := -\tfrac{1}{\rho} \ln \left(\sum\textstyle_{j=1}^{N_c} e^{-\rho h_j(q, \omega)}\right),\label{eq:cbf_example_electromechanical}
\end{align}
where $\rho > 0$ modulates the approximation sharpness. 
The reduced-order safety controller is then designed to satisfy the robust CBF condition
\begin{align}
  \nabla h(q, \omega)\T f_{\text{rs}}(q, \omega) + \gamma h(q, \omega) - \rob \geq 0,
\end{align}
where $f_{\text{rs}}$ is the vector field of the reduced-order system~\eqref{eq:motor_pendulum_red}. 
We implement a smooth safety controller, following the approach discussed in \cite{cohen2023characterizing}, which is applied to a nominal PD controller with gravity compensation. The system is initialized at rest ($q(0) = [0^\circ, 0^\circ]^\top$, $\omega(0) = [0, 0]^\top$) and is controlled toward the target configuration $q_\text{target} = [180^\circ, 0^\circ]^\top$. 
The safety controller parameters are $\gamma = 10$, $\rob = 0.1$, and $\rho = 50$.
Simulations are performed for $10$ linearly spaced values of $\pr \in \left[0.001, 0.035\right]$.
\begin{figure}[htpb]
    \centering
    \begin{minipage}[c]{0.49\columnwidth}
        \centering
        \includegraphics[width=\linewidth]{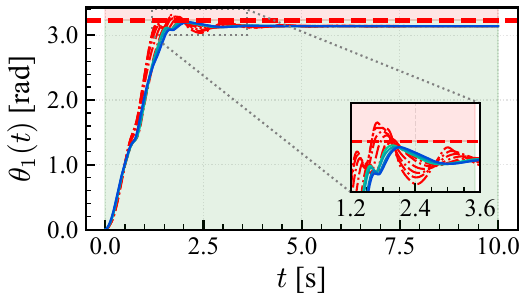}
        \label{fig:theta1_dp}
    \end{minipage}
    \hfill
    \begin{minipage}[c]{0.49\columnwidth}
        \centering
        \includegraphics[width=\linewidth]{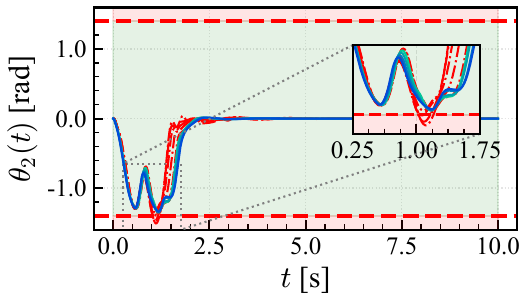}
        \label{fig:theta2_dp}
    \end{minipage}

    \vspace{-1em}

    \begin{minipage}[c]{0.49\columnwidth}
        \centering
        \includegraphics[width=\linewidth]{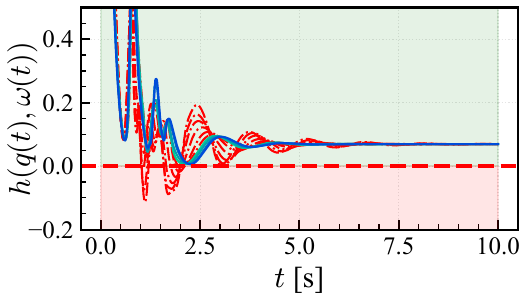}
        \label{fig:h_dp}
    \end{minipage}
    \hfill
    \begin{minipage}[c]{0.49\columnwidth}
        \centering
        \includegraphics[width=\linewidth]{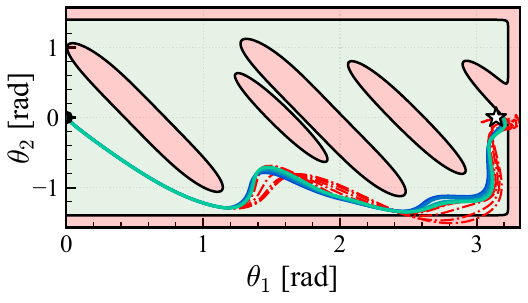}
        \label{fig:cspace_dp}
    \end{minipage}
    \vspace{-1em}
    \caption{Results for the 2-DoF robotic arm with joint motors. Blue-to-cyan solid lines indicate safe trajectories
    ($\pr<\bar{\pr}$); dash-dotted red
    lines denote safety violations ($\pr\geq\bar{\pr}$).
    \textbf{Top left:} Time evolution of $\theta_1$.
    \textbf{Top right:} Time evolution of $\theta_2$.
    \textbf{Bottom left:} Time evolution of the CBF \eqref{eq:cbf_example_electromechanical}.
    \textbf{Bottom right:} Corresponding trajectories in the joint space. The red areas represent forbidden configurations induced by the obstacles in the robot workspace.}
    \label{fig:dp_results}
\end{figure}
Results are reported in Fig. \ref{fig:dp_results}: specifically, the top left panel and the top right panel illustrate the corresponding joint trajectories, the bottom left panel shows the evolution of the aggregated CBF \eqref{eq:cbf_example_electromechanical}, and the bottom right panel depicts the system trajectories in the joint configuration space, highlighting the regions associated with the prescribed safety constraints. 
For sufficiently small values of $\pr$, the trajectories remain within the safe region, consistent with the timescale separation paradigm established in Theorem~\ref{th:SP}.
By contrast, increasing $\pr$ slows the electrical dynamics relative to the mechanical subsystem, leading to larger transient discrepancies between the actual motor currents and their equilibrium values.
As illustrated, these mismatches may lead to safety violations, despite the controller being designed to enforce robust safety on the reduced-order dynamics.

\subsection{Physical Dynamics with Safety Control Algorithm}

  We now consider the suboptimal safe control scenario introduced in Section \ref{sec:dynamic_safe_control}, where timescale separation arises from the interaction between a physical plant and a dynamic optimization algorithm.
  Unlike the previous example, the fast subsystem does not represent neglected physical dynamics, but rather the dynamic refinement of the control input through a suitable optimization algorithm that asymptotically imposes the safety constraint.
  The virtual experiments are performed on a 6-DoF UR10e robotic manipulator simulated in the MuJoCo physics engine.
  The robot is modeled at the kinematic level, with joint configuration $q \in \R^6$ and joint velocity $\dot{q} = u \in \R^6$.
  Starting from the initial configuration $q(0) = [0, -\frac{\pi}{2}, \frac{\pi}{2}, -\frac{\pi}{2}, -\frac{\pi}{2}, 0]^\top$, the manipulator is tasked with following a joint-space reference path $q_\text{ref}(s)$, obtained through piecewise cubic interpolation of a set of waypoints.
  The nominal control policy is defined as
  \begin{align*}
  \pi^{\text{des}}(q, s) &= K_p (q_\text{ref}(s) - q) + \frac{\partial q_\text{ref}(s)}{\partial s} \dot{s},
  \end{align*}
  with $K_p = 10 I_6$. The path-progress variable $s \in [0, 1]$ is governed by the dynamics
  
  \begin{equation}
  \dot{s} =
  \begin{cases}
  \displaystyle
  \frac{v_s}
  {1+\|q_{\mathrm{ref}}(s)-q\|^2/e_s^2},
  & 0\leq s<1,\\[6pt]
  0, & s=1,
  \end{cases}
\end{equation}
  with parameters $v_s = 0.4$ and $e_s = 5.0$.

\begin{figure}[htpb]
    \centering
    \begin{minipage}[c]{0.35\columnwidth}
        \centering
        \includegraphics[width=\linewidth]{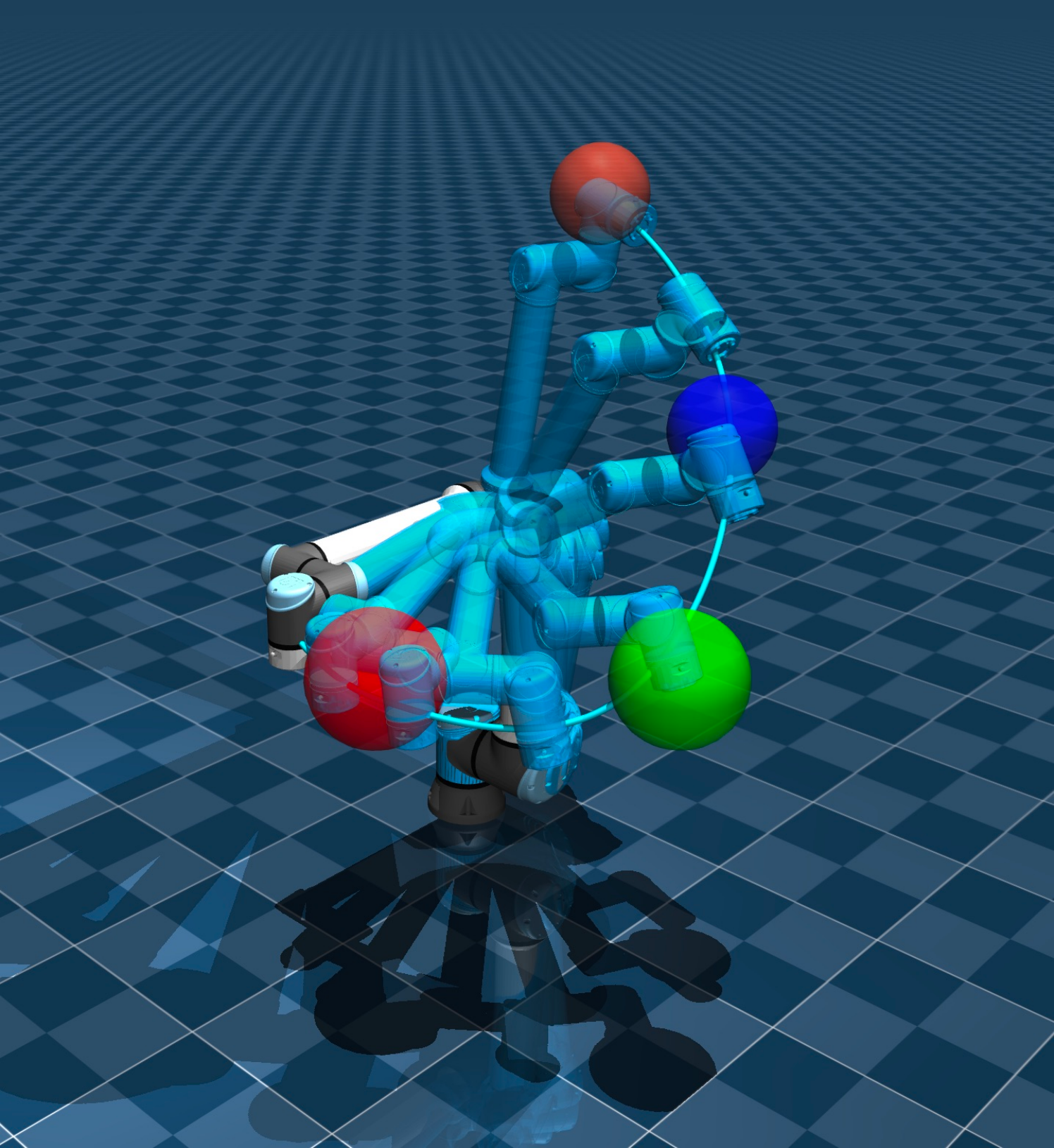}
        \label{fig:nominal_ee_ur}
    \end{minipage}%
    % \hfill%
    \hspace{0.05\columnwidth}%
    \begin{minipage}[c]{0.35\columnwidth}
        \centering
        \includegraphics[width=\linewidth]{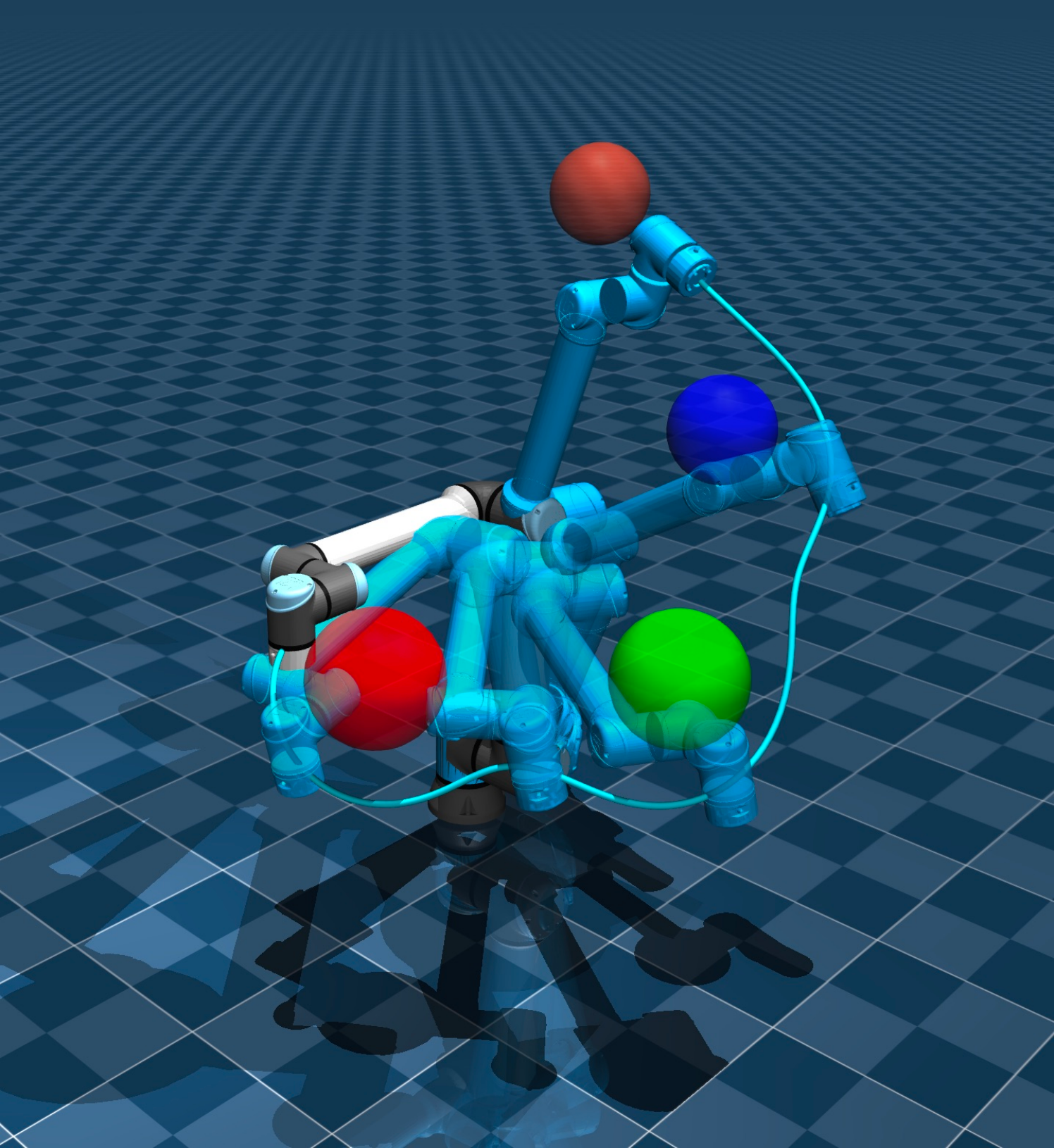}
        \label{fig:nominal_trail_ur}
    \end{minipage}
    \caption{\textbf{Left:} Nominal unsafe trajectory. \textbf{Right:} Safe trajectory with the suboptimal safety filter ($\pr = 0.01$).
    }
    \label{fig:nominal_ur}
\end{figure}

The environment contains $N_o = 4$ spherical obstacles, deliberately positioned near the nominal end-effector path to induce collision-avoidance maneuvers. See the left panel of Fig. \ref{fig:nominal_ur} for the nominal (i.e., without safety constraints) end-effector trajectory and selected robot configurations along the nominal path.
To account for possible collisions involving the tool flange of the robot, we represent its position by a point $p(q) \in \R^{3}$ and account for its physical dimensions by introducing an inflation radius $r_{\mathrm{p}} = 0.085$ m.
For each obstacle $j \in \until{N_o}$, we define the safety function
\begin{align*} 
  \bar{h}_{j}(q) := \norm{p(q) - c_j}^2 - (r_j + r_{\mathrm {p}})^2,
\end{align*} 
where $c_j \in \R^{3}$ and $r_j > 0$ are the center and radius of the $j$-th spherical obstacle, respectively.
Moreover, we incorporate the joint limits in the safety function. Given the lower and upper joint limits $q_i^{\min}, q_i^{\max} \in \R$, we define
\begin{align}
  \chi_{\cQ}(q) := \prod_{i=1}^{6}\psi\left(
    \frac{2q_i-q_i^{\max}-q_i^{\min}}{q_i^{\max}-q_i^{\min}}
  \right),
\end{align}
where $\psi: \R \to [0,1]$ is a smooth function defined as
\begin{align}
  \psi(c)=
  \begin{cases}
    \exp\left(-\frac{c^2}{1-c^2}\right)
    & |c|<1
    \\
    0 &|c|\geq1.
  \end{cases}
\end{align}
The resulting safety function $h: \R^6 \to \R$ is then defined as
\begin{align}
  h(q) := \beta \left( \chi_{\mathcal Q}(q) \prod_{j=1}^{N_o} \bar{h}_j(q) - \tau \right), \label{eq:cbf_ur10e}
\end{align}
where $\beta, \tau > 0$ are tunable design parameters.
Specifically, $\beta$ is a scaling factor, and $\tau$ selects a compact safe set with a regular boundary.
Using the safety function \eqref{eq:cbf_ur10e}, the safety filtering problem is formulated as 
\begin{align}
  \min_{u \in \R^6} &\quad \norm{u - \pi^{\text{des}}(q, s)}^2 \notag\\
  \text{s.t.} &\quad \nabla h(q)^\top u + \gamma h(q) - \rob \geq 0.\label{eq:cbf_constraint_ur10e}
\end{align}
In accordance with the proposed framework, both the joint configuration $q$ and the path-progress variable $s$ constitute the slow state $x:=\mathrm{col}(q,s)\in\mathbb R^6\times[0,1]$. 
Rather than solving \eqref{eq:cbf_constraint_ur10e} to optimality at every instant, we dynamically adjust the control input through a gradient-flow algorithm applied to an unconstrained version of problem~\eqref{eq:cbf_constraint_ur10e}.
Specifically, we introduce a smooth cost function $J: \R^6 \times [0,1] \times \R^6 \to \R$ defined as
\begin{align}
  J(x, u) &:= \frac{1}{2}\norm{u - \pi^{\text{des}}(x)}^2 + \kappa(x) S_\delta(-\res(x, u)), \label{eq:cost_function_ur10e}
\end{align}
where %
$\res(x, u) := \nabla h(q)^\top u + \gamma h(q) - \rob$ encodes the robust safety constraint and $S_\delta: \R \to \R$ is the smooth Softplus function
\begin{align}
  S_\delta(\bar{c}) &:= \delta \log(1 + e^{\frac{\bar{c}}{\delta}}), \label{eq:softplus}
\end{align}
with $\delta > 0$ a tunable parameter that controls the approximation sharpness.
The state-dependent weight $\kappa: \R^6 \times [0,1] \to \R$ is
\begin{align}
  \kappa(x):= \frac{2\left(S_\delta(-\resd(x))+\nu\right)} {\norm{\nabla h(q)}^2+\vartheta(\resd(x))}, \label{eq:state_dependent_weight}
\end{align}
where $\resd(x) := \nabla h(q)^\top \pi^{\text{des}}(x) + \gamma h(q) - \rob$ is the residual of the robust CBF condition obtained with the nominal control input, $\nu > 0$ is a tuning parameter, and $\vartheta: \R \to \R$ reads as
\begin{align}
  \vartheta(\bar{c}):=
  \begin{cases}
  0,&\bar{c}\leq0
  \\
  c_0\exp\left(-(r_0/\bar{c})^2\right),&\bar{c}>0.
  \end{cases}
\end{align}
Thus, by introducing $\pr$ to modulate the learning rate, the resulting gradient-flow dynamics are
\begin{align}
  \pr \dot{u} &= - (u - \pi^{\text{des}}(x)) - \kappa(x) \nabla_u S_\delta(-c(x, u)).
  \label{eq:grad_flow_ur10e}
\end{align}
Hence, in this scenario, the control input constitutes the fast state, i.e., $z := u \in \R^6$.
In particular, since $J(x,u)$ is strongly convex in $u$ uniformly in $x$, system~\eqref{eq:grad_flow_ur10e} admits an exponentially stable equilibrium manifold $\{(x,u)| u = u^\star(x)\}$, where $u^\star(x) := \arg\min_{u \in \R^6} J(x, u)$, thus satisfying Assumption~\ref{ass:bl}.
Hence, in this scenario, the reduced-order model corresponds to the plant dynamics controlled with $u = u^\star(x)$.
Then, since it is possible to show that, for sufficiently large values of $\beta$ and $\tau$, the equilibrium $u^\star(x)$ strictly satisfies the robust CBF condition, i.e., $c(x, u^\star(x)) > 0$ for all $x \in \{(q,s) \in \R^{6} \times [0,1] \mid \cbf(q) \geq 0\}$, Assumption~\ref{ass:cbf} is also satisfied.
The interconnected plant and gradient-flow dynamics are numerically integrated using the Radau method available in SciPy, with the solver restarted at each outer simulation time step $\Delta t = 0.005$ s. 
The parameters are set to $\gamma = 10$, $\rob = 0.01$, $\nu = 0.01$, $c_0 = 1$, $r_0 = 0.1$, $\delta = 0.005$, $\beta=1$, $\tau=10^{-6}$.
Simulations are performed over a $5$-second horizon for 19 values of $\pr \in [0.01, 1]$, chosen as 9 linearly spaced values in $[0.01, 0.09]$ and 10 linearly spaced values of $\pr \in [0.1, 1]$.
To evaluate the effect of finite optimization speed, we also consider an ideal safety controller obtained by solving the optimization problem~ \eqref{eq:cbf_constraint_ur10e} to optimality at each time step, using IPOPT within the CasADi framework~\cite{andersson2019casadi}. 
\begin{figure}[H]
  \centering
  \includegraphics[width=0.48\textwidth]{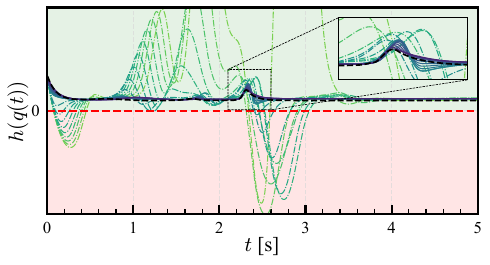}
  \caption{Time evolution of the CBF \eqref{eq:cbf_ur10e} for the scenario involving physical dynamics and a safety control algorithm. The viridis color gradient indicates increasing values of $\pr$. Safe trajectories are represented by solid lines, whereas unsafe trajectories are represented by dash-dotted lines.
  The optimal solution is represented by the dashed black line.
  }
  \label{fig:h_evolution}
\end{figure}
Figure \ref{fig:h_evolution} compares the time evolution of the aggregated safety function~\eqref{eq:cbf_ur10e} for the selected values of $\pr$. 
For sufficiently small values of $\pr$, the gradient-flow algorithm evolves on a faster timescale than the manipulator, allowing the computed suboptimal control input to track its instantaneous equilibrium more closely. See, e.g., the resulting trajectory for $\pr = 0.01$ in the right panel of Fig. \ref{fig:nominal_ur}.
Correspondingly, the resulting trajectories satisfy the prescribed safety constraints. 
As $\pr$ increases, the optimization algorithm responds more slowly to changes in the robot configuration, producing larger transient deviations from its equilibrium. 
The numerical results show that these deviations can lead to violations of the safety constraints.

\section{Conclusions}

  In this paper, we investigated safety of nonlinear feedback-interconnected systems through a singular perturbation approach. 
  By suitably combining reduced-order safety properties and stability conditions of the fast dynamics, we established sufficient conditions under which a controller designed solely on the reduced-order model guarantees safety of the overall interconnection.
  Furthermore, we characterized the dependence of the certified safe sets on the timescale separation parameter, establishing their nesting property and asymptotic recovery of the reduced-order safe region as the fast dynamics become increasingly fast.
  We illustrated the applicability of the proposed framework through numerical simulations on an electromechanical 2-DoF robotic arm and a suboptimal safe control architecture in which a 6-DoF robotic manipulator interacts with an online optimization algorithm.

\appendix

\subsection{Proof of Theorem~\ref{th:SP}. Part I: Forward Invariance} 
\label{sec:proof2}

The idea of the proof is to combine the reduced-system safety properties encoded in $\cbf$ (cf. Assumption~\ref{ass:cbf}) with the stability properties of the boundary-layer system~\eqref{eq:bl} encoded in $U$ (cf. Assumption~\ref{ass:bl}) in order to guarantee that $V_\pr(x,z)$ ensures safety for the overall interconnected system~\eqref{eq:interconnected_system_generic}.
To this end, the proof exploits the parameter $\pr$ to accelerate the fast dynamics relative to the slow dynamics.
Intuitively, this allows one to make the approximation errors arbitrarily small. These errors arise from the fact that the properties of $\cbf$ and $U$ are established only for the limit cases detailed in Assumptions~\ref{ass:cbf} and~\ref{ass:bl}, rather than for the original interconnection~\eqref{eq:interconnected_system_generic}.

First of all, since they will be useful later, let us introduce the function $\zeqcl: \cC \to \R^{p}$ and the constants $\nablahbound, \fredbound, \lip_{\zeqcl}  > 0$ defined as
\begin{subequations}
  \begin{align}
    \zeqcl(x) &:= \zeq(x,\pol(x))
    \label{eq:zeqcl}
    \\
    \lip_{\zeqcl} &:= \max\left\{\delta, \sup_{x \in \cC}\norm{\nabla \zeqcl(x)}\right\}\label{eq:zeqcl_lip}
    \\
    \fredbound &:= \sup_{x \in \cC}  
    \norm{\sdyn(\x,\zeqcl(\x),\pol(\x))} .
    \label{eq:bound_phi}
    \\
    \nablahbound &:= \sup_{x \in \cC}  
    \norm{\nabla \cbf(x)}
    \label{eq:nabla_W_alpha_W}
\end{align}
\end{subequations}
We highlight that the above constants exist and are finite since $\cC$ is compact (cf. Assumption~\ref{ass:cbf}) and $\zeq$, $\nabla\zeq$, $\nabla \cbf$, $\pol$, and $\sdyn$ are continuous (see Assumptions~\ref{ass:equilibrium},~\ref{ass:cbf}, and~\ref{ass:lip}, respectively).
Moreover, the arbitrary constant $\delta>0$ ensures that $L_{\zeqcl} > 0$, even when the equilibrium manifold is constant, thereby avoiding degeneracies in the parameter bounds introduced below.
Then, by Lipschitz continuity of $\sdyn$ and $\fdyn$ (cf. Assumption~\ref{ass:lip}) and $\pol$ (cf. Assumption~\ref{ass:cbf}), we note that the solution to system~\eqref{eq:interconnected_system_generic} exists and is unique (see, e.g.,~\cite[Thm.~3.2]{khalil2002nonlinear}).
Let us rewrite system~\eqref{eq:interconnected_system_generic} in the error coordinate $\tzt := \zt - \zeq(\xt,\pol(\xt))$, thus obtaining 
\begin{subequations}\label{eq:interconnected_system_generic_err}
  \begin{align}
    \xtp &= \scl(\xt,\tzt)
    \label{eq:slow_system_generic_err}
    \\
    \tztp &= \frac{1}{\pr}\fcl(\tzt,\xt) -\nabla\zeqcl(\xt)\xtp
    ,\label{eq:fast_system_generic_err}
  \end{align}
\end{subequations}
in which we introduce $\scl: \cC \times \R^{p} \to \R^{n}$ and $\fcl: \R^p \times \cC \to \R^{p}$ defined as 
\begin{subequations}
  \begin{align}
   \scl(x,\tz) &:= \sdyn(x,\tz + \zeqcl(x),\pol(x))
    \label{eq:scl}
    \\
    \fcl(\tz,x) &:= \fdyn(\tz + \zeqcl(x),x,\pol(x)).
    \label{eq:fcl}
  \end{align}
\end{subequations}
Then, we evaluate the Lie derivative of $\cbf$ along the trajectories of subsystem~\eqref{eq:slow_system_generic_err}.
Specifically, by adding and subtracting the term $\nabla\cbf(\x)\T\drs(\x,\pol(\x)) = \nabla\cbf(\x)\T\scl(\x,0)$ (see the definitions of $\drs$~\eqref{eq:reduced} and $\scl$~\eqref{eq:scl}), we obtain
\begin{align}
  \nabla\cbf(\x)\T\scl(\x,\tz)
  &= 
  \nabla\cbf(\x)\T\drs(\x,\pol(\x))
  % \notag\\
  % &\hspace{.4cm}
  + \nabla\cbf(\x)\T(\scl(\x,\tz)-\scl(\x,0))
  \notag\\
  &\stackrel{(a)}{\geq}
  -\alpha(\cbf(\x)) + \rob
  % \notag\\
  % &\hspace{.4cm}
  - \norm{\nabla\cbf(\x)}\norm{\scl(\x,\tz)-\scl(\x,0)},
  \label{eq:deltaW_1}
\end{align}
for all $(x,z) \in \cC \times \R^p$, where in $(a)$ we exploit~\eqref{eq:CBF} (cf. Assumption~\ref{ass:cbf}) and we use the Cauchy-Schwarz inequality.
Since the vector field $f$ is $\lip_{\sdyn}$-Lipschitz continuous (cf. Assumption~\ref{ass:lip}), the difference $\scl(\x,\tz)-\scl(\x,0)$ satisfies 
\begin{align}\label{eq:lipschitz_bound}
  \norm{\scl(\x,\tz)-\scl(\x,0)} \leq \lip_{\sdyn} \norm{\tz}.
\end{align} 
By using inequality~\eqref{eq:lipschitz_bound} and $\lip_{\cbf}$ defined in~\eqref{eq:nabla_W_alpha_W}, we are able to further bound the right-hand side of~\eqref{eq:deltaW_1} as
\begin{align}
  \nabla\cbf(\x)\T\scl(\x,\tz)
  &\geq -\alpha(\cbf(\x)) + \rob - \lip_{\cbf}\lip_{\sdyn} \norm{\tz}.\label{eq:deltaW}
\end{align}
Now, we evaluate the Lie derivative of the Lyapunov function $U$ characterized in Assumption~\ref{ass:bl} along the trajectories of~\eqref{eq:fast_system_generic_err} and we use~\eqref{eq:U_minus_generic} to write
\begin{align}
  \frac{1}{\pr}\nabla U(\tz)\T\fcl(\tz,\x) - \nabla U(\tz)\T\nabla\zeqcl(\x)\xtp
  % \notag\\
  &
  \leq 
  - \frac{b_3}{\pr}\norm{\tz}^2 - \nabla U(\tz)\T\nabla\zeqcl(\x)\xtp
  \notag\\        
  &\stackrel{(a)}{\leq}
  - \frac{b_3}{\pr}\norm{\tz}^2 + b_4\lip_{\zeqcl}\norm{\tz}\norm{\scl(\x,\tz)}
  \notag\\        
  &\stackrel{(b)}{\leq}
  - \frac{b_3}{\pr}\norm{\tz}^2 + b_4\lip_{\zeqcl}\norm{\tz}\norm{\drs(\x,\pol(\x))}
  % \notag\\
  % &\hspace{.4cm}
  +b_4\lip_{\zeqcl}\norm{\tz}\norm{\scl(\x,\tz) - \scl(\x,0)}\notag\\        
  &\stackrel{(c)}{\leq}
  - \frac{b_3}{\pr}\norm{\tz}^2 + b_4\lip_{\zeqcl}\fredbound\norm{\tz}
  +b_4\lip_{\zeqcl}\lip_{\sdyn}\norm{\tz}^2,
  \label{eq:DeltaU}
\end{align} 
where in $(a)$ we use the Cauchy-Schwarz inequality, the bounds~\eqref{eq:lipshcitz_hessian_u} (cf. Assumption~\ref{ass:bl}) and~\eqref{eq:zeqcl_lip}, and~\eqref{eq:slow_system_generic_err}, in $(b)$ we add and subtract the term $\drs(\x,\pol(\x)) = \scl(x,0)$ (see the definitions of $\drs$~\eqref{eq:reduced} and $\scl$~\eqref{eq:scl}) in the last norm, and use the triangle inequality, while in $(c)$ we use the Lipschitz continuity of $\sdyn$ (cf. Assumption~\ref{ass:lip}) and the bound~\eqref{eq:bound_phi}.
Now, let us define some constants that will be useful in the following, namely
\begin{equation}
    K_1 = \nablahbound \lip_{\sdyn}, \quad
    K_2 = b_4 \lip_{\zeqcl} \lip_{\sdyn}, \quad
    K_3 = b_4 \lip_{\zeqcl} \fredbound.
\end{equation}
Consequently, we can rewrite the bounds~\eqref{eq:deltaW} and~\eqref{eq:DeltaU} as
\begin{subequations}
  \begin{align}
    \dot{\cbf}(\x)
    &\geq -\alpha(\cbf(\x)) + \rob - K_1 \norm{\tz} \label{eq:deltaW_2}
    \\
    \dot{U}(\tz)
    &\leq - \left(\frac{b_3}{\pr} - K_2\right) \norm{\tz}^2 + K_3 \norm{\tz}.\label{eq:DeltaU_2}
  \end{align}
\end{subequations}
We now introduce the scalar function $r: \cC \times \R^p \to \R$ 
\begin{align}
  r(x, z) := \frac{K_1}{\sqrt{b_1}}\sqrt{U(z - \zeqcl(x))}.\label{eq:r}
\end{align}
For brevity, we will denote $r(\tz) := r(x, \tz + \zeqcl(x))$ in the following.
By using the lower bound in~\eqref{eq:U_first_bound_generic} (cf. Assumption~\ref{ass:bl}) and the definition of $r$ (cf.~\eqref{eq:r}), we have
\begin{equation}
  K_1 \norm{\tz} \leq \frac{K_1}{\sqrt{b_1}}\sqrt{U(\tz)} = r(x,z).\label{eq:r_bound}
\end{equation}
Consequently, by substituting \eqref{eq:r_bound} into \eqref{eq:deltaW_2}, we obtain
\begin{align}
  \dot{\cbf}(\x) &\geq -\alpha(\cbf(\x)) + \rob - r(x,z).\label{eq:deltaW_3}
\end{align}
Now, we characterize the evolution of $r(x,z)$. 
To this end, let us define $\bar{\pr}_1 := \frac{b_3}{K_2}$.
Then, for all $\pr \in (0,\bar{\pr}_1)$, we use the bound~\eqref{eq:U_first_bound_generic} to further bound the right-hand side of~\eqref{eq:DeltaU_2} as
\begin{align}
  \dot{U}(\tz) \leq -\frac{1}{b_2}\left(\frac{b_3}{\pr} - K_2\right) U(\tz) + \frac{K_3}{\sqrt{b_1}} \sqrt{U(\tz)}.
  \label{eq:DeltaU_3}
\end{align}
Moreover, since $U(\tz) \geq 0$ (see the lower bound~\eqref{eq:U_first_bound_generic} in Assumption~\ref{ass:bl}), we have
\begin{align}
  \frac{d}{dt} \sqrt{U(\tz)} = \frac{\dot{U}(\tz)}{2\sqrt{U(\tz)}},
\end{align}
for all $\tz \in \R^p$ such that $U(\tz) > 0$.
Therefore, by substituting~\eqref{eq:DeltaU_3} into the above equation, we obtain
\begin{align}
  \frac{d}{dt} \sqrt{U(\tz)} \leq - \frac{1}{2b_2}\left(\frac{b_3}{\pr} - K_2\right) \sqrt{U(\tz)} + \frac{K_3}{2\sqrt{b_1}}.\label{eq:DeltaU_4}
\end{align}
By multiplying both sides by $\frac{K_1}{\sqrt{b_1}}$ and using the definition of $r(x,z)$ (cf.~\eqref{eq:r}) in~\eqref{eq:DeltaU_4}, we obtain 
\begin{align}
  \dot{r}(x,z) &\leq - \frac{K_1}{2b_2\sqrt{b_1}}\left(\frac{b_3}{\pr} - K_2\right) \sqrt{U(\tz)} + \frac{K_1 K_3}{2 b_1}
  \notag\\
  &= - \frac{1}{2b_2}\left(\frac{b_3}{\pr} - K_2\right) r(x,z) + \frac{K_1 K_3}{2 b_1}.
  \label{eq:dr_bound}
\end{align}
Now, let us introduce 
\begin{align}
  \lp &:= \frac{1}{2b_2}\left(\frac{b_3}{\pr} - K_2\right)
  \label{eq:lp_def}
  \\
  d &:= \frac{K_1 K_3}{2 b_1}.\label{eq:d_def}
\end{align}
These definitions allow us to compactly rewrite~\eqref{eq:dr_bound} as
\begin{align}
  \dot{r}(x,z) &\leq - \lp r(x,z) + d.
  \label{eq:dr_bound_2}
\end{align}
Then, we introduce the upper bound $\bar{\pr}$ defined as
\begin{align}
  \bar{\pr} := \frac{b_1 b_3 \rob}{b_2 K_1 K_3 + b_1 K_2 \rob}.
\end{align}
First of all, we note that $\bar{\pr} < \bar{\pr}_1$. Now, let us consider an arbitrary $\pr \in (0, \bar{\pr})$. Consequently, the bound~\eqref{eq:dr_bound_2} holds and, by defining $\rho_{\pr} := \frac{d}{\lp}$, it holds
\begin{align}\label{eq:rho_rob}
  \rho_{\pr} < \rob,
\end{align}
for all $\pr \in (0, \bar{\pr})$.
With this result at hand, we introduce $\xipr: \R_{\ge 0} \to \R_{\ge 0}$ defined as
\begin{align}
  \xipr(r) \!:=\! 
  \begin{cases}
    0 \!\!& \text{if } r \!\leq\! \rob
    \\
    \frac{(r - \rho_{\pr})^{\lip_{\alpha}/\lp}}{\lp} \int_{\rob}^{r} (s - \rob)(s - \rho_{\pr})^{-\lip_{\alpha}/ \lp - 1} ds \!\!& \text{if } r \!>\! \rob,
  \end{cases}
  \label{eq:xi}
\end{align}
and use it to construct the composite barrier function $V_\pr: \cC \times \R^p \to \R$ defined as
\begin{align}
  V_\pr(x, z) := \cbf(x) - \xipr(r(z - \zeqcl(x))).\label{eq:V_pr}
\end{align}
To lighten the notation, we also define 
$\tilde{V}_\pr: \cC \times \R^p \to \R$ as
\begin{align}
  \tilde{V}_\pr(x, \tz) := V_\pr(x, \tz + \zeqcl(x)),\label{eq:tilde_V}
\end{align}
as it will be useful in the following to analyze the evolution of $V_\pr$ along the trajectories of~\eqref{eq:interconnected_system_generic_err}.
From the definition of $\xipr$ (cf.~\eqref{eq:xi}), we obtain by direct computation
\begin{align}
  \xipr^\prime(r) = \begin{cases}
    0 & \text{if } r \leq \rob
    \\
    \frac{L_\alpha \xipr(r) + r - \rob}{\lp r - d} & \text{if } r > \rob,
  \end{cases}\label{eq:xi_derivative_def}
\end{align}
which allows us to write
\begin{align}\label{eq:derivative_xi}
  (\lp r - d) \xipr^\prime(r) - \lip_{\alpha} \xipr(r) = r - \rob,
\end{align}
for all $r > \rob$.
Since $U$ is continuously differentiable and positive definite (cf. Assumption~\ref{ass:bl}), $r$ is continuously differentiable wherever $r > 0$. 
Moreover, by construction, $\xipr(r) = 0$ for $r \leq \rob$ and $\xipr(r) = \xipr^\prime(r) = 0$ for $r = \rob$. Therefore, since $\cbf$ is also continuously differentiable, $\tilde{V}_\pr$ is continuously differentiable for all $(x, \tz) \in \cC \times \R^p$ such that $\tilde{V}_\pr(x, \tz) \geq 0$. 
Now, consider an arbitrary pair $(x, \tz)$ such that $\tilde{V}_\pr(x, \tz) \geq 0$.
Since $\xipr(r) \geq 0$, we have that $\cbf(x) \geq 0$, which implies $x \in \cC$. 
We now distinguish two cases: (i) $r(x,z) \leq \rob$ and (ii) $r(x,z) > \rob$. 
Suppose that case (i) holds and, thus, $\xipr(r(\tz)) = 0$ and $\xipr^\prime(r(\tz)) = 0$. 
Therefore, by using~\eqref{eq:deltaW_3}, in this first case, the Lie derivative of $\tilde{V}_\pr$ along the trajectories of~\eqref{eq:interconnected_system_generic_err} satisfies
\begin{align}
  \dot{\tilde{V}}_\pr(x, \tz) 
  &\geq -\alpha(\cbf(x)) + \rob - r
  % \notag\\
  \stackrel{(a)}{=}  -\alpha(\tilde{V}_\pr(x, \tz)) + \rob - r
  % \notag\\
  \stackrel{(b)}{\geq} -\alpha(\tilde{V}_\pr(x, \tz)),
  \label{eq:lie_V_generic}
\end{align}
where in $(a)$ we use the definition of $\tilde{V}_\pr$ (cf.~\eqref{eq:tilde_V}) and the fact that $r(x,z) \leq \rob \implies \xipr(r) = 0$, while in $(b)$ we use the fact that $r(x,z) \leq \rob$.
Consider now case (ii) and, thus, $\xipr(r(x,z)), \xipr^\prime(r(x,z)) > 0$.
By using~\eqref{eq:deltaW_3} and~\eqref{eq:dr_bound_2}, we have
\begin{align}
  \dot{\tilde{V}}_\pr(x, \tz)
  &\geq -\alpha(\cbf(x)) + \rob - r + \xipr^\prime(r)(\lp r - d).\label{eq:Lie_V}
\end{align}
Since $\alpha$ is $\lip_{\alpha}$-Lipschitz continuous and increasing (cf. Assumption~\ref{ass:cbf}) and $\xipr(r(x,z)) \geq 0$ (see~\eqref{eq:xi}), it holds
\begin{align*}
  \alpha(\cbf(x)) \leq \alpha(\cbf(x) - \xipr(r)) + \lip_{\alpha} \xipr(r),
\end{align*}
which, combined with the definition of $\tilde{V}_\pr$ (cf.~\eqref{eq:tilde_V}), allows us to further bound the right-hand side of~\eqref{eq:Lie_V} as
\begin{align}
  \dot{\tilde{V}}_\pr(x, \tz)
  &\geq -\alpha(\tilde{V}_\pr(x, \tz)) - \lip_{\alpha} \xipr(r) + \rob - r 
  % \notag\\
  % &\hspace{.4cm}
  + \xipr^\prime(r) (\lp r - d) 
  % \notag\\
  \stackrel{(a)}{=} -\alpha(\tilde{V}_\pr(x, \tz)),\label{eq:last}
\end{align}
where in $(a)$ we use~\eqref{eq:derivative_xi}.
  Finally, we verify the nonvanishing gradient condition required by Lemma \ref{lemma:cbf}. Consider an arbitrary $(x,z) \in \partial \cCf(\pr)$. 
  If $r(x,z) \leq \rob$, then $V_\pr(x,z) = \cbf(x) = 0$.
  By the strict safety margin in Assumption \ref{ass:cbf}, we have $\nabla \cbf(x) \neq 0$, implying that $\nabla V_\pr(x,z) \neq 0$.
  If instead $r(x,z) > \rob$, then $\tz \neq 0$, and Assumption \ref{ass:bl} guarantees that $\nabla U(\tz) \neq 0$.
  Moreover, $\xipr^\prime(r) > 0$ follows from \eqref{eq:xi_derivative_def} and the fact that $\rho_\pr < \rob$. Hence,
  \begin{align*}
    \nabla_2 V_\pr(x,z) = - \frac{K_1 \xipr^\prime(r)}{2 \sqrt{b_1 U(\tz)}} \nabla U(\tz) \neq 0.
  \end{align*}
  Therefore, we conclude that $\nabla V_\pr(x,z) \neq 0$ for all $(x,z) \in \partial \cCf(\pr)$.
The proof follows by using~\eqref{eq:last} to invoke Lemma~\ref{lemma:cbf}.

\subsection{Proof of Theorem~\ref{th:SP}. Part II: Nesting of Certified Sets}\label{sec:proof_cor1}

The first part of the proof consists in showing that
\begin{align}
  \xi_{\pr_1}(r) \leq \xi_{\pr_2}(r),\quad \forall r \geq 0, \quad \forall \pr_1 \leq \pr_2.\label{eq:xi_monotone}
\end{align}
To this end, by rearranging~\eqref{eq:derivative_xi} (cf. Appendix~\ref{sec:proof2}), we write
\begin{align}\label{eq:xi_prime}
  \xipr^\prime(r) = \frac{\lip_{\alpha} \xipr(r) + r - \rob}{\lp r - d},
\end{align}
for all $r > \rob$.
Moreover, from the definition of $\xipr$ in \eqref{eq:xi}, we recall that $\xipr(r) = 0$ for all $r \leq \rob$.
Now, let us consider $0 < \pr_1 \leq \pr_2 < \bar{\pr}$ and, to lighten the notation, let us define $\lambda_1 := \lambda_{\pr_1}$ and $\lambda_2 := \lambda_{\pr_2}$, where $\lambda_\pr$ for a generic $\pr$ is defined in~\eqref{eq:lp_def} (cf. Appendix~\ref{sec:proof2}).
This definition implies
\begin{align}\label{eq:lambda_monotone}
  \pr_1 \leq \pr_2 \implies \lambda_1 \geq \lambda_2.
\end{align}
Furthermore, given $\rho_{\pr_i} := \frac{d}{\lambda_i}$, the bound~\eqref{eq:rho_rob} leads to
\begin{align}
  \lambda_i r - d > 0, \quad \forall r \geq \rob.
\end{align} 
For $r \leq \rob$, the desired inequality \eqref{eq:xi_monotone} holds because 
\begin{align*}
  \xi_{\pr_1}(r) = \xi_{\pr_2}(r) = 0.
\end{align*}
For $r > \rob$, let us define 
\begin{align}\label{eq:Delta}
  \Delta(r) := \xi_{\pr_2}(r) - \xi_{\pr_1}(r).
\end{align}
To achieve~\eqref{eq:xi_monotone}, we need to prove $\Delta(r) \geq 0$ for all $r > \rob$.
Differentiating $\Delta(r)$ and using \eqref{eq:xi_prime}, we obtain
\begin{align}
  \Delta^\prime(r) = \xi_{\pr_2}^\prime(r) - \xi_{\pr_1}^\prime(r)
  % \notag\\
  &= \frac{\lip_{\alpha} \xi_{\pr_2}(r) + r - \rob}{\lambda_2 r - d} - \frac{\lip_{\alpha} \xi_{\pr_1}(r) + r - \rob}{\lambda_1 r - d}\notag\\
  &\stackrel{(a)}{=} \frac{\lip_\alpha}{\lambda_2 r - d} \Delta(r)
  % \notag\\
  % &\hspace{.4cm}
   + \left(
    \frac{1}{\lambda_2 r - d} - \frac{1}{\lambda_1 r - d}
  \right) \left(
    \lip_\alpha \xi_{\pr_1}(r) + r - \rob
  \right),\label{eq:Delta_prime}
\end{align}
where in $(a)$ we add and subtract the term $\frac{\lip_\alpha \xi_{\pr_1}(r) + r - \rob}{\lambda_2 r - d}$.
Now, we recall that $\lambda_1 \geq \lambda_2$ (see~\eqref{eq:lambda_monotone}); thus, the second term on the right-hand side of the above equation is non-negative. Observe also that $\xi_{\pr_1}(r) \geq 0$ by construction and that $r - \rob > 0$. 
Therefore, by using~\eqref{eq:Delta_prime}, we get the bound
\begin{align}
  \Delta^\prime(r)  - \frac{\lip_\alpha}{\lambda_2 r - d} \Delta(r) \geq 0.\label{eq:delta_prime}
\end{align}
Then, let us introduce $\mu: [\rob, \infty) \to \R_+$ defined as
\begin{align}
  \mu(r) := \exp\left(
    - \int_{\rob}^{r} \frac{\lip_\alpha}{\lambda_2 s - d}ds
  \right) > 0,\label{eq:mu}
\end{align}
which is strictly positive by positivity of the exponential function and by $\lambda_2 r - d > 0$ for all $r \geq \rob$.
Thus, it holds
\begin{align}
  \frac{d}{dr}\mu(r) \Delta(r) &= 
  \mu(r)\left(-\frac{\lip_\alpha}{\lambda_2 r - d}\right) \Delta(r) + \mu(r)\Delta^\prime(r)
  \stackrel{(a)}{\geq} 0,
  \label{eq:delta_prime_2}
\end{align}
where in $(a)$ we use~\eqref{eq:delta_prime} and $\mu(r) \ge 0$.
Moreover, we observe
\begin{align}\label{eq:integral}
 \int_{\rob}^{r} \frac{d}{ds}\mu(s) \Delta(s) ds =
  \mu(r) \Delta(r) - \mu(\rob) \Delta(\rob) \geq 0.
\end{align}
Since $\Delta(\rob) = 0$ (see the definitions of $\Delta$ and $\xipr$ in~\eqref{eq:Delta} and~\eqref{eq:xi}, respectively) and $\mu(r) > 0$ for all $r \geq \rob$ (see~\eqref{eq:mu}), \eqref{eq:integral} implies $\Delta(r) \geq 0$ for all $r > \rob$, thus proving \eqref{eq:xi_monotone}. 
Consequently, recalling the definition of $V_\pr$ (cf.~\eqref{eq:V_pr}), we have
\begin{align}
  V_{\pr_2}(x,z) &= 
  h(x) - \xi_{\pr_2}(r(z - \zeq(x,\pol(x)))) 
  \\
  &\leq h(x) - \xi_{\pr_1}(r(z - \zeq(x,\pol(x)))) = V_{\pr_1}(x,z).\label{eq:last_cor}
\end{align}
The proof follows by combining~\eqref{eq:last_cor} and the definition of the certified safe set $\cCf(\pr)$ (cf.~\eqref{eq:certified_safe_set}).

\subsection{Proof of Theorem~\ref{th:SP}. Part III: Asymptotic Recovery}\label{sec:proof_cor2}

First, let us establish that $V_\pr$ converges to $\cbf$ uniformly on every compact subset of $\cC \times \R^p$ as $\pr \to 0^+$.
Let $\cK \subseteq \cC \times \R^p$ be an arbitrary nonempty compact set.
By definition of $V_\pr$ (cf.~\eqref{eq:V_pr}) and since $\xipr$ is nonnegative (see its definition in~\eqref{eq:xi}), we have
\begin{align}\label{eq:first}
   \cbf(x) - V_\pr(x,z) = \xipr(r(z - \zeq(x,\pol(x)))) \geq 0.
\end{align}
Since $r$ and $\zeq(x, \pi(x))$ are continuous functions (see the definition of $r$ in~\eqref{eq:r} and Assumptions~\ref{ass:equilibrium},~\ref{ass:cbf}, and~\ref{ass:bl}), there exists a positive constant $R_\cK < \infty$ such that $r(z - \zeq(x, \pi(x))) \leq R_\cK$ for all $(x,z) \in \cK$.
If $R_\cK \leq \rob$, the desired equality holds by definition of $\xipr$ (cf.~\eqref{eq:xi}).
Otherwise, if $R_\cK > \rob$, let us define $M_\cK := R_\cK - \rob > 0$. 
Then, by using~\eqref{eq:derivative_xi}, we have
\begin{align}
  \xipr^\prime(r) = \frac{\lip_{\alpha} \xipr(r) + r - \rob}{\lp r - d} \leq \frac{\lip_{\alpha} \xipr(r) + r - \rob}{\lp \rob - d},\label{eq:xi_prime_bound}
\end{align}
for all $r \in [\rob, R_\cK]$ and $\pr \in (0, \bar{\pr})$. 
Since $\xipr \geq 0$ and $\lp \rob - d > 0$ for $\pr \in (0, \bar{\pr})$, the bound~\eqref{eq:xi_prime_bound} implies $\xipr^\prime(r) \geq 0$, thus proving that $\xipr$ is non-decreasing in $r$.
Integrating both sides of \eqref{eq:xi_prime_bound} over $[\rob, r]$ (with $r \in [\rob, R_\cK]$) and using $\xipr(\rob) = 0$ (see~\eqref{eq:xi}), we obtain
\begin{align}
  \xipr(r) ^= \int_{\rob}^{r} \xipr^\prime(s) ds
  &\stackrel{(a)}{\leq}
  \frac{1}{D_\pr} \int_{\rob}^{r} (\lip_{\alpha} \xipr(s) + s - \rob) ds 
  \notag\\
  &\stackrel{(b)}{\leq}
  \frac{1}{D_\pr} \int_{\rob}^{r} (\lip_{\alpha} \xipr(r) + s - \rob) ds 
  \notag\\
  &\stackrel{(c)}{=} \frac{1}{D_\pr} \left[ \lip_{\alpha} \xipr(r) (r - \rob) + \frac{(r - \rob)^2}{2} \right]
  \notag\\
  &\stackrel{(d)}{\leq}
   \frac{L_\alpha M_\cK \xipr(r) + M_\cK^2/2}{D_\pr},\label{eq:xi_bound}
\end{align}
where in $(a)$ we use the bound~\eqref{eq:xi_prime_bound} and introduce $D_\pr:= \lp \rob - d$, in $(b)$ we use the fact that $\xipr$ is non-decreasing, in $(c)$ we evaluate the integral, while in $(d)$ we use the fact that $r - \rob \leq M_\cK$.
By recalling the definition of $\lp$ in \eqref{eq:lp_def}, we have
\begin{align}
  D_\pr = \frac{\rob b_3}{2b_2 \pr} - \frac{K_2 \rob}{2b_2} - d.
  \label{eq:D_pr}
\end{align}
By rearranging~\eqref{eq:xi_bound} and using~\eqref{eq:D_pr}, we obtain
\begin{align}
  \left(
    \frac{\rob b_3}{2b_2 \pr} - \frac{K_2 \rob}{2b_2} - d - L_\alpha M_\cK
  \right) \xipr(r) \leq \frac{M_\cK^2}{2}.\label{eq:bound}
\end{align}
Given $\bar{\pr}_{\cK,1} := \min\{\bar{\pr}, \frac{\rob b_3}{4b_2 (\frac{K_2 \rob}{2b_2} + d + L_\alpha M_\cK)}\}$, the bounds~\eqref{eq:xi_bound} and~\eqref{eq:bound} ensure that, for all $\pr \in (0, \bar{\pr}_{\cK,1})$, it holds
\begin{align}\label{eq:cck_pr}
  0 \leq \xipr(r) \leq \frac{M_\cK^2}{2(D_\pr - L_\alpha M_\cK)} \leq c_\cK \pr,
\end{align}
where $c_\cK := \frac{2 b_2 M_\cK^2}{\rob b_3} > 0$. 
Combining~\eqref{eq:cck_pr} and~\eqref{eq:first} yields
\begin{align}
  0 \leq h(x) - V_\pr(x,z) \leq c_\cK \pr,\label{eq:uniform_convergence}
\end{align}
for all $(x,z) \in \cK$ and $\pr \in (0, \bar{\pr}_{\cK,1})$. 
Taking the supremum over $\cK$ in \eqref{eq:uniform_convergence} and letting $\pr \to 0^+$ establishes the claimed uniform convergence, that is,
\begin{align}
  \lim_{\pr \to 0^+} \sup_{(x,z) \in \cK} | V_\pr(x,z) - \cbf(x) | = 0.
\end{align}
Now, let $\cortwoset$ be an arbitrary compact subset of $\{(x,z) \in \cC \times \R^p \mid h(x) > 0\}$.
As $\cortwoset$ is a particular example of the compact set $\cK$ considered above, we can apply the same arguments used above to find $c_{\cortwoset} > 0$ and $\bar{\pr}_{\cortwoset,1} > 0$ such that
\begin{align}
  0 \leq h(x) - V_\pr(x,z) \leq c_{\cortwoset} \pr,\label{eq:uniform_convergence_2}
\end{align}
for all $(x,z) \in \cortwoset$ and $\pr \in (0, \bar{\pr}_{\cortwoset,1})$.
Moreover, by continuity of $h$ and compactness of $\cortwoset$, there exists $\delta_\cortwoset := \min_{(x,z) \in \cortwoset} h(x) > 0$. From \eqref{eq:uniform_convergence_2}, we thus obtain
\begin{align}
  V_\pr(x,z) %
  \geq \delta_\cortwoset - c_\cortwoset \pr.
\end{align}
Therefore, given $\bar{\pr}_{\cortwoset} := \min\{\bar{\pr}_{\cortwoset,1}, \frac{\delta_\cortwoset}{c_{\cortwoset}}\}$, we have $V_\pr(x,z) > 0$ for all $\pr \in (0, \bar{\pr}_{\cortwoset})$ and $(x,z) \in \cortwoset$.
This implies that $\cortwoset \subseteq \cCf(\pr)$. %
Hence, we are ready to prove the recovery of the entire safe region.
Since $\cCf(\pr) \subseteq \cC \times \R^p$ for all $\pr \in (0, \bar{\pr})$ and since $\cC \times \R^p$ is closed (cf. Assumption~\ref{ass:cbf}), we have
\begin{align}
  \overline{
    \bigcup_{\pr \in (0, \bar{\pr})}\cCf(\pr) 
  } \subseteq \cC \times \R^p. \label{eq:closure1}
\end{align}
To establish the reverse inclusion, consider any $(x,z) \in \cC \times \R^p$ such that $\cbf(x) > 0$.
Applying the previous arguments to the singleton $\cortwoset = \{(x,z)\}$ yields $(x,z) \in \cCf(\pr)$ for sufficiently small values of $\pr$.
If instead $\cbf(x) = 0$, Assumption~\ref{ass:cbf} yields
\begin{align}
  \nabla \cbf(x)^\top \drs(x, \pol(x)) \geq \rob > 0.
\end{align}
Setting $v:= \drs(x, \pol(x))$ and given $s > 0$, the continuous differentiability of $\cbf$ (cf. Assumption~\ref{ass:cbf}) allows us to perform a first-order Taylor expansion of $\cbf(x + sv)$ around $x$ and obtain
\begin{align}
  \cbf(x + sv) =s \nabla \cbf(x)^\top v + o(s) > 0,
\end{align}
for sufficiently small values of $s$. 
With these results at hand, let us consider a sequence of points $\{(x + s_k v, z)\}_{k \in \N}$ with $s_k > 0$ and $s_k \to 0^+$ as $k \to \infty$.
By applying the previous argument to each singleton $\cortwoset_k = \{(x + s_k v, z)\}$, we show that every point $(x + s_k v, z)$ belongs to $\cCf(\pr)$ for sufficiently small values of $\pr > 0$.
Since $(x + s_k v, z) \to (x,z)$ as $k \to \infty$, we then use the definition of closure to obtain
\begin{align}
  (x,z) \in \overline{
    \bigcup_{\pr \in (0, \bar{\pr})}\cCf(\pr) 
  }.
\end{align}
Since $(x,z) \in \cC \times \R^p$ was arbitrary, we conclude that every point in $\cC \times \R^p$ belongs to the closure of the union of the sets $\cCf(\pr)$, i.e.,
\begin{align}
  \cC \times \R^p \subseteq \overline{
    \bigcup_{\pr \in (0, \bar{\pr})}\cCf(\pr) 
  }.\label{eq:closure2}
\end{align}
The combination of~\eqref{eq:closure1} and~\eqref{eq:closure2} yields the desired result.

\end{document}